\documentclass[a4paper,11pt]{article}

\usepackage{authblk}
\usepackage{amsmath}
\usepackage{amsfonts}
\usepackage{amssymb}
\usepackage{latexsym}
\usepackage{indentfirst}
\usepackage[top=1in, bottom=1in, left=1in, right=1in]{geometry}
\usepackage{algorithm}
\usepackage{algpseudocode}
\usepackage{graphicx,psfrag,epsf}
\usepackage{xcolor}
\usepackage[font=small]{caption}
\usepackage{subcaption} 
\usepackage{float}    
\usepackage{enumerate}
\usepackage{enumitem}
\usepackage{natbib}
\usepackage{url}
\usepackage{bm}
\usepackage{longtable}
\usepackage{multirow}

\let\footnote=\endnote

\newtheorem{theorem}{Theorem}

\newtheorem{remark}{Remark}
\newtheorem{corollary}{Corollary}

\allowdisplaybreaks

\def\T{{ \mathrm{\scriptscriptstyle T} }}

\DeclareMathOperator*{\EE}{E}
\DeclareMathOperator*{\mse}{MSE}
\DeclareMathOperator{\imse}{IMSE}

\DeclareMathOperator*{\hsgp}{HSGP}

\DeclareMathOperator*{\diag}{diag}

\DeclareMathOperator*{\dist}{dist}

\DeclareMathOperator*{\tr}{tr}

\def\om{\mathbf{\om}}
\def\si{\mathbf{\si}}

\newcommand{\ud}{\mathrm{d}}

\def\C{\mathbf{C}}

\def\G{\mathbf{G}}

\def\I{\mathbf{I}}

\def\K{\mathbf{K}}

\def\W{\mathbf{W}}
\def\X{\mathbf{X}}
\def\Y{\mathbf{Y}}

\DeclareMathOperator{\Lcal}{\mathcal{L}}
\DeclareMathOperator{\Mcal}{\mathcal{M}}

\def\RR{\mathbb{R}}
\def\NN{\mathbb{N}}
\def\ZZ{\mathbb{Z}}

\def\EE{\mathbb{E}}
\def\PP{\mathbb{P}}

\def\bfc{\mathbf{c}}

\def\bfh{\mathbf{h}}

\def\bfj{\mathbf{j}}
\def\bfk{\mathbf{k}}

\def\bft{\mathbf{t}}

\def\bfw{\mathbf{w}}
\def\bfx{\mathbf{x}}

\def\bfz{\mathbf{z}}

\begin{document}

\title{\bf Fast and Provably Accurate Sequential Designs using Hilbert Space Gaussian Processes}

\author[1]{Huanyan Zhu}
\author[1]{Cheng Li}
\affil[1]{Department of Statistics and Data Science, National University of Singapore}
\date{}
\maketitle

\begin{abstract}
Gaussian processes are widely used for accurate emulation of unknown surfaces in sequential design of expensive simulation experiments. Integrated mean squared error (IMSE) is an effective acquisition function for sequential designs based on Gaussian processes. However, existing approaches struggle with its implementation because the required integrals often lack closed-form expressions for most kernel functions. We propose a novel and computationally efficient Hilbert space Gaussian process approximation for the IMSE acquisition function, where a truncated eigenbasis representation of the integral enables closed-form evaluation. We establish sharp global non-asymptotic bounds for both the approximation error of isotropic kernels and the resulting error in the acquisition function. In a series of numerical experiments with $\gamma$-stabilizing, the proposed method achieves substantially lower prediction error and reduced computation time compared to existing benchmarks. These results demonstrate that the proposed Hilbert space Gaussian process framework provides an accurate and computationally efficient approach for Gaussian process based sequential design.
\end{abstract}


\section{Introduction} \label{sec:intro}
Sequential design with Gaussian process (GP) surrogates has become an important and increasingly popular approach for learning expensive response surfaces under limited evaluation budgets. Rather than pre-specifying all sampling locations, the method selects new inputs iteratively based on the current surrogate fit, yielding a principled trade-off between exploration and exploitation. This approach is particularly useful when each function evaluation is costly, time-consuming, or otherwise difficult to obtain, as in ocean-acoustic source localization (\citealt{jenkins:2023}), hybrid studies involving both computer models and field data (\citealt{koermer:2024}), and additive manufacturing process optimization (\citealt{ding:2025}).

Consider sequential design over a bounded domain $\Omega\subset\RR^d$. After observing data $\{(\bfx_n, y_n)\}_{n=1}^N$ with $\bfx_n\in\Omega, y_n \in \RR$, Gaussian process regression provides a flexible surrogate with accurate prediction and principled uncertainty quantification. At each step, a new input is selected by optimizing an acquisition function $a_N(\bfx)$:
\begin{align}
  \bfx_{N+1} \in \arg\min_{\bfx \in \Omega} a_N(\bfx), \nonumber 
\end{align}
where $a_N$ is constructed from the GP predictive distribution and tailored to the experimental objective. Broadly speaking, existing criteria fall into two categories: those that aim to improve global predictive accuracy, such as \emph{integrated mean squared error/integrated mean squared prediction error} (IMSE/IMSPE) (\citealt{binois:2019}, \citealt{koermer:2024}, \citealt{ding:2025}), and those that aim to optimize the unknown response surface, such as \emph{expected improvement} (EI) (\citealt{jenkins:2023}; \citealt{chen:2023}) and \emph{entropy search} (ES) (\citealt{weichert:2024}). The latter typically allocate evaluations near extrema and do not explicitly enforce predictive accuracy across the whole domain. In contrast, IMSE/IMSPE criteria target global error reduction by integrating the posterior mean squared error over $\Omega$, which promotes more uniform refinement of the surrogate throughout the design space. For this reason, IMSE is especially appealing in applications that require accurate emulation over the full domain, including sensitivity analysis and uncertainty propagation. Accordingly, the focus of this study is on IMSE-based sequential design.

We now introduce the IMSE-based sequential design criterion. At each iteration, the next input is chosen to maximize the expected reduction in integrated posterior uncertainty:
\begin{align} \label{eq:intro_imse}
   \bfx_{N+1} \in \arg\max_{\bfx \in \Omega} \Big\{{\imse}_N - {\imse}_{N+1}(\bfx) \Big\},
\end{align}
In this expression, $\imse_N$ denotes the current integrated posterior mean squared error, while $\imse_{N+1}(\bfx)$ denotes the IMSE obtained after augmenting the design with a hypothetical observation at $\bfx$. This one-step lookahead criterion directly quantifies the global utility of a candidate point by measuring its effect on reducing predictive uncertainty over $\Omega$. IMSE-based designs have a long history in the literature on computer experiments and surrogate modeling, where they are valued for their strong global predictive performance.

Despite its effectiveness in reducing global prediction error, evaluating and optimizing \eqref{eq:intro_imse} is often computationally demanding. The criterion involves an integral over $\Omega$, and the integrand depends on the candidate location $\bfx$ through the GP update. Consequently, naive evaluation by dense quadrature or Monte Carlo integration within an outer optimization loop can be prohibitively expensive. On the other hand, closed-form or highly optimized expressions are available in restricted cases, typically for special combinations of kernels and integration measures.
For example, IMSPE-based sequential design has been studied and implemented efficiently for several widely used kernels in \citet{binois:2019}. Their approach rewrites the IMSPE criterion so that the candidate dependence is absorbed into matrix updates, leaving the computation of quantities of the form $W_{ij} := \int k(\bfx_i, \bfx) k(\bfx, \bfx_j) \ud \PP(\bfx)$ where $k$ is a kernel function. Although closed-form expressions for $W_{ij}$ are available for Gaussian, Mat\'ern-$3/2$ and Mat\'ern-$5/2$ kernels with uniform measure $\PP$, extending such calculus-based or kernel-specific derivations to broader families can be nontrivial. As a result, practical implementations often remain restricted to a small number of covariance functions. Accordingly, the \texttt{hetGP} package \citep{binois:2021} provides efficient IMSPE computation for sequential designs, but its implementation also relies on precomputed kernel integrals and currently supports only the Gaussian, Mat\'ern-$3/2$, Mat\'ern-$5/2$ kernels. 
Recent work has extended IMSPE-type criteria beyond the classical computer-experiment setting, but mainly through additional problem-specific analytic derivations rather than through a general resolution of the underlying computational bottleneck. For example, \citet{koermer:2024} develop an IMSPE criterion for hybrid computer-model and field-data experiments, deriving closed-form expressions and gradients that enable sequential augmentation in that specialized setting. Similarly, \citet{ding:2025} propose an active learning stochastic kriging procedure for additive manufacturing reliability and stability analysis, using an IMSPE-style criterion to guide sequential sampling.
While these developments broaden the applicability of IMSPE in important domains, they do not resolve the central computational challenge. In these cases, efficient IMSPE evaluation still relies on kernel/measure-specific analytic calculations, which restrict practical implementations to a small number of covariance families and make generalization to broader kernels difficult.

The intractable kernel integral problem also appears in kernel mean embeddings and Bayesian quadrature (BQ). In particular, standard kernel mean embedding concerns the first-order quantity $\int k(\bfx_i, \bfx) \ud \PP(\bfx)$. While this form differs from the second-order kernel integrals required for the IMSE, which takes the form $\int k(\bfx_i, \bfx) k(\bfx, \bfx_j) \ud \PP(\bfx)$, existing research in mean embeddings provides a valuable foundation for addressing the IMSE problem. For instance, Table 1 in \citet{briol:2019} summarizes a collection of tractable kernel means in certain conjugate settings, which extends the limited set of Gaussian, Mat\'ern-3/2, Mat\'ern-5/2 kernels in \citet{binois:2021}. Despite these advances, the general case remains intractable and therefore motivates approximation. Several lines of work have sought either to enlarge the tractable class or to approximate non-analytic embeddings. For shift-invariant kernels, \citet{nishiyama:2016} derives convolution-based representations for both kernel mean values and reproducing kernel Hilbert space (RKHS) inner products, and proposes convolution-trick constructions that reduce these quantities to the evaluation of densities within an infinitely divisible family. However, tractability is not guaranteed in general, since many infinitely divisible densities do not admit closed-form expressions and may still require numerical evaluation. Alternatively, when samples from $\PP$ are available, sampling-based approaches approximate the mean embedding through low-rank representations, such as Nystr\"om-type estimators with carefully chosen representative points and weights (\citealt{chatalic:2022}, \citealt{hayakawa:2023}). In another direction, \citet{warren:2024} develops fast Fourier transform based approximations for stationary kernels, combining Fourier-domain transforms with interpolation operators to approximate kernel means efficiently.
Nevertheless, these developments do not directly resolve the IMSE integral arising in sequential design. First, the IMSE integrand involves a kernel product, which is inherently more restrictive for analytic solutions than single-kernel embeddings. Second, purely sampling-based estimators can be prohibitively expensive in sequential design, because the IMSE integrand changes at every iteration, making repeated Monte Carlo or Nystr\"om approximations computationally burdensome. Consequently, there remains a critical need for an approximation strategy that balances high accuracy with the computational efficiency required for iterative design.

To address this computational bottleneck, this work proposes a fast and accurate approximation to IMSE based on a \emph{Hilbert space Gaussian process} (HSGP) representation. HSGP methods approximate Gaussian processes through truncated eigenfunction expansions of the Laplace operator on a bounded domain, with coefficients determined by the kernel spectral density (\citealt{solin:2020}). In practice, HSGP has been shown to provide accurate and computationally efficient approximations to kernel functions (\citealt{riutort:2023}). This structure yields an efficient surrogate for the IMSE objective that avoids expensive numerical integration while maintaining high fidelity to the original criterion. In particular, the proposed HSGP approach applies to a broad class of kernels with closed-form spectral densities, including the Gaussian kernel, generalized Wendland kernels, and Mat\'ern kernels beyond a few half-integer smoothness cases, and offers substantial computational savings in sequential design settings where the acquisition function must be evaluated repeatedly. In addition, detailed convergence rate analysis is established for the HSGP approximation of both the kernel and the $\imse$ criterion.
\vspace{2mm}

The main contributions are summarized as follows:
\vspace{-2mm}
\begin{itemize}[leftmargin=5mm]
  \item \textbf{Method:} An HSGP-based acquisition function $\widehat{\imse}_m(\bft)$ is developed by combining a truncated Laplacian eigenbasis expansion with an analytic reduction of the IMSE integral to a quadratic form, where $m$ denotes the number of eigenfunctions. In contrast to the IMSPE criterion in \citet{binois:2019}, which supports only the Gaussian, Mat\'ern-$3/2$, and Mat\'ern-$5/2$ kernels, the proposed HSGP approach applies to any kernel family with a closed-form spectral density and avoids numerical integration in sequential design.
  
  \item \textbf{Theory:} Explicit approximation bounds are derived for the underlying kernel approximation and then propagated to acquisition-level error guarantees. A $\gamma$-stabilizing feasibility restriction is introduced to ensure numerical stability by controlling the power-function denominator, and this restriction also implies quasi-uniformity of the resulting sequential designs. While \citet{solin:2020} establish polynomial-decay error bounds for HSGP kernel approximation with implicit constants, our work strengthens these results by deriving exponential-order error bounds with explicit constants. In addition, a new error bound is established for the HSGP-based approximation of $\imse$-type acquisition functions.
  
  \item \textbf{Empirical evaluation:} One and two dimensional experiments compare the HSGP-IMSE strategy with the IMSPE baseline implemented in the \texttt{hetGP} package when available, and with a space-filling baseline based on Latin hypercube sampling otherwise. The comparison examines predictive accuracy, uncertainty contraction, and runtime over sequential iterations. The numerical results show that the HSGP-IMSE approach outperforms the competing methods in reducing both the global \emph{root mean square prediction error} (RMSE) and the global mean posterior predictive variance, while achieving a runtime comparable to that of IMSPE. These findings demonstrate that the proposed HSGP-IMSE method is both accurate and computationally efficient.
\end{itemize}

The remainder of this paper is organized as follows. Section~\ref{sec:sqgp} reviews GP regression and IMSE-based sequential design, introduces the HSGP approximation, and presents the resulting computational procedure for $\widehat{\imse}_m(\bft)$ together with a $\gamma$-stabilizing strategy. Section~\ref{sec:theory} develops theoretical guarantees for the approximation errors of both the HSGP kernel representation and the IMSE acquisition function. Section~\ref{sec:experiments} presents numerical experiments that evaluate approximation fidelity as well as end-to-end sequential design performance.

\section{Sequential design using Gaussian process} \label{sec:sqgp}
This section presents a standard framework for sequential design with GP regression. Section~\ref{subsec:sequential_gp} reviews GP models and introduces the $\imse$ acquisition function, together with maximum likelihood estimation of the GP hyperparameters. Section~\ref{subsec:HSGP-IMSE} introduces an $\hsgp$-based approximation to the IMSE acquisition function, leading to a closed-form approximation of $\imse$ and making its optimization computationally tractable. Section~\ref{subsec:gamma_greedy} discusses a $\gamma$-stabilizing strategy for robust sequential sampling.

\subsection{Gaussian process regression} \label{subsec:sequential_gp}
Let the design domain be $\Omega=(-B,B)^d \subset \RR^d, B > 0$. Given inputs $\X_N := \{\bfx_1,\ldots,\bfx_N\}\subset \Omega$ and observations $\Y_N := [y_1,\ldots,y_N]^\T \in \RR^N$, consider the GP regression model
\begin{align} \label{eq:model_gpr}
    y_n &= \Mcal(\bfx_n) + \epsilon_n, \quad n = 1,\ldots,N, \\
    \Mcal(\cdot) &\sim \mathcal{GP}(\mu(\cdot), k(\bfx, \bfx^\prime)), \nonumber\\
    \epsilon_n &\overset{\text{i.i.d.}}{\sim} \mathcal{N}(0,\eta), \qquad \eta \ge 0, \nonumber
\end{align}
where $\mathcal{GP}(\mu(\cdot), k(\bfx, \bfx^\prime))$ is a Gaussian process with mean function $\mu(\cdot)$ and covariance function $k:\Omega\times\Omega\to\RR$, where $k$ is a positive definite kernel, and $\{\epsilon_n\}_{n=1}^N$ are independent of $\Mcal$. Here $\eta=0$ corresponds to the noiseless/interpolatory setting, while $\eta>0$ allows either Gaussian observation noise or a fixed nugget regularization term. We assume that $\mathbf{X}_N$ contains only distinct points. If replicate inputs are present, they should be averaged. 
Model \eqref{eq:model_gpr} is treated as our working GP model, while the response may in fact be generated from a mechanism of the form $y=f(\bfx)+\epsilon$, where the underlying function $f$ is unknown. The main objective of sequential design is to construct an accurate GP surrogate for $f$ over the entire domain $\Omega$, thereby providing reliable emulation for all $\bfx\in\Omega$ (\citealt{binois:2019}).

For shift-invariant kernels, the covariance function depends only on the difference between two inputs, that is, $k(\bfx,\bfx') = k(\bfx-\bfx')$. Two standard examples are the Gaussian kernel and the Mat\'ern kernel.
\begin{itemize}
    \item The squared exponential (Gaussian) kernel with length scale $\ell$ and scale $\sigma^2$,
    \begin{align} \label{eq:kernel_gaussian}
        k(\bfx, \bfx^\prime) = \sigma^2 \exp\left(-\frac{\|\bfx-\bfx^\prime\|_2^2}{2 \ell^2}\right).
    \end{align}
    \item The Mat\'ern kernel with length scale $\ell$, scale $\sigma^2$, and smoothness $\nu > 0$,
    \begin{align} \label{eq:kernel_matern}
        k(\bfx, \bfx^\prime) = \sigma^2 \frac{2^{1-\nu}}{\Gamma(\nu)} \left(\frac{\sqrt{2\nu} \|\bfx-\bfx^\prime\|_2}{\ell}\right)^\nu K_\nu\left(\frac{\sqrt{2\nu} \|\bfx-\bfx^\prime \|_2}{\ell}\right),
    \end{align}
    where $K_\nu$ is the modified Bessel function of the second kind.
\end{itemize}

\paragraph{Hyperparameter estimation.} Hyperparameter estimation is carried out using a modeling strategy similar to that of \citet{binois:2019}. Assume that the mean function $\mu(\bfx)$ is specified and treated as known.
Define $\Y_N = \left[y_1, \ldots, y_N \right]^\T$, $\bm\mu_N = [\mu(\bfx_1),\ldots,\mu(\bfx_N)]^\T$, $\K_N := (k(\bfx_i,\bfx_j))_{1\le i,j\le N}$, and write the kernel as $k(\cdot,\cdot)=\sigma^2 c_\ell(\cdot,\cdot)$ with correlation function $c_\ell$ and correlation matrix $\C_N:=\bigl(c_\ell(\bfx_i,\bfx_j)\bigr)_{ij}$. Similarly, the nugget is reparameterized as $\eta=\sigma^2 g$, so that $\K_N + \eta \I_N=\sigma^2(\C_N + g\I_N)$. With this reparameterization, the Gaussian log-likelihood is
\begin{align} \label{eq:log-likelihood}
    \Lcal(\sigma^2,\ell, g)
    &= \mathrm{Const}
    -\frac{N}{2}\log(\sigma^2)
    -\frac{1}{2}\log|\C_N + g \I_N| \nonumber\\
    &\qquad -\frac{1}{2\sigma^2}\bigl(\Y_N-\bm\mu_N\bigr)^\T \left(\C_N + g \I_N \right)^{-1} \bigl(\Y_N-\bm\mu_N\bigr).
\end{align}
For fixed $(\ell,g)$, the Maximum likelihood estimation (MLE) of $\sigma^2$ admits closed form:
\begin{align} \label{eq:mle_sigma2}
    \hat \sigma^2
    &=
    \frac{1}{N}\bigl(\Y_N-\bm\mu_N\bigr)^\T \left(\C_N + g \I_N \right)^{-1} \bigl(\Y_N-\bm\mu_N\bigr).
\end{align}
By contrast, the maximum likelihood estimators of $\ell$ and $g$ do not admit closed-form expressions and must be obtained numerically, for example by gradient-based optimization. The gradients with respect to $\ell$ and $g$ are
\begin{align} \label{eq:mle_numerical}
    \frac{\partial}{\partial \ell} \Lcal
    &=
    -\frac{1}{2}\tr \left[\left(\C_N + g \I_N \right)^{-1}\frac{\partial \C_N}{\partial \ell}\right] \nonumber\\
    &\qquad+\frac{1}{2\sigma^2} \bigl(\Y_N-\bm\mu_N\bigr)^\T
    \left(\C_N + g \I_N \right)^{-1}\frac{\partial \C_N}{\partial \ell} \left(\C_N + g \I_N \right)^{-1}
    \bigl(\Y_N-\bm\mu_N\bigr), \nonumber\\
    \frac{\partial}{\partial g} \Lcal
    &=
    -\frac{1}{2}\tr \left[\left(\C_N + g \I_N \right)^{-1}\right] + \frac{1}{2\sigma^2} \bigl(\Y_N-\bm\mu_N\bigr)^\T
    \left(\C_N + g \I_N \right)^{-2} \bigl(\Y_N-\bm\mu_N\bigr). 
\end{align}
In implementations, $\sigma^2$ in \eqref{eq:mle_numerical} will be replaced by \eqref{eq:mle_sigma2} during optimization over $\ell$ and $g$.

\paragraph{IMSE acquisition.} Define $\X_N := \{\bfx_n\}_{n=1}^N$, $\bfk_N(\bfx) := [k(\bfx,\bfx_1),\ldots,k(\bfx,\bfx_N)]^\T$ and $\K_N := (k(\bfx_i,\bfx_j))_{1\le i,j\le N}$. Suppose the mean function $\mu(\bfx)$ is known. 
From model \eqref{eq:model_gpr}, the posterior mean predictor of latent process $\mathcal{M}(\bfx)$ at an arbitrary $\bfx\in \Omega$ is
\begin{align} 
    \widehat{\mathcal{M}}(\bfx) = \mu(\bfx) + \bfk_N^\T(\bfx) \left(\K_N + \eta \I_N \right)^{-1} \left( \Y_N-\bm\mu_N \right). \nonumber
\end{align}
The posterior mean squared error (MSE) for predicting the latent process $\mathcal{M}(\bfx)$ is
\begin{align} 
    \mse(\bfx \mid \X_N)
    &= \EE \left[(\mathcal{M}(\bfx) - \widehat{\mathcal{M}}(\bfx))^2 \mid \X_N, \Y_N \right]
      = k(\bfx,\bfx) - \bfk_N^\T(\bfx) \left(\K_N + \eta \I_N \right)^{-1} \bfk_N(\bfx). \nonumber
\end{align}
Let $\PP$ be a finite measure on $\Omega$ such as the Lebesgue measure or the uniform probability measure. For a candidate point $\bft\in\Omega$, define the reduction in integrated MSE as 
\begin{align} 
    \imse(\bft)
    := \int_{\Omega}\left[\mse(\bfx\mid \X_N)-\mse(\bfx\mid \X_N\cup\{\bft\})\right] \ud \PP(\bfx). \nonumber
\end{align}
The following identity expresses $\imse(\bft)$ as an integrated squared residual scaled by a power function.

\begin{theorem} \label{thm:imse_acq}
Assume $\PP$ is a finite measure on $\Omega$ and define the power function
\begin{align} \label{eq:power_function}
    P_{N, \eta}(\bft)
    :=\left[k(\bft,\bft)-\bfk_N^\T(\bft) \left(\K_N + \eta \I_N \right)^{-1} \bfk_N(\bft)\right]^{1/2}.
\end{align}
Suppose $P_{N, \eta}^2(\bft) + \eta > 0$, then
\begin{align} \label{eq:imse_int}
    \imse(\bft)
    = \frac{1}{P_{N, \eta}^2(\bft) + \eta}
    \int_{\Omega}\left[k(\bfx,\bft)-\bfk_N^\T(\bfx) \left(\K_N + \eta \I_N \right)^{-1}\bfk_N(\bft)\right]^2  \ud \PP(\bfx).
\end{align}
\end{theorem}

The IMSE-type acquisition admits a transparent decomposition that highlights its built-in balance between \emph{exploration} and \emph{exploitation}. On one hand, the denominator $P_{N, \eta}^2(\bft)$ normalizes the global impact term and prevents the criterion from being dominated purely by local uncertainty, which encourages \emph{exploration}. On the other hand, the numerator quantifies the integrated discrepancy between the prior correlation pattern $k(\bft,\cdot)$ and its projection onto the span of the observed kernel sections $\{k(\bfx_n, \cdot)\}_{n=1}^N$. Large numerator thus favors candidates that are expected to yield a substantial reduction of the overall predictive uncertainty, which can be viewed as \emph{exploitation} of the current surrogate by refining regions that most impact the global prediction quality.

\paragraph{Sequential sampling algorithm.} A baseline sequential design procedure based on the exact acquisition $\imse(\bft)$ is summarized in Algorithm~\ref{alg:sequential_with_exact_imse}. The main computational bottleneck is the evaluation of $\imse(\bft)$ over candidates $\bft\in\Omega$, which typically requires numerical integration (e.g., Monte Carlo or quadrature) together with repeated linear algebra involving $\K_N$. In sequential sampling, this evaluation must be carried out at every iteration, and the overall cost quickly becomes prohibitive as the number of samples $N$ increases. This motivates the development of fast and accurate surrogates for $\imse(\bft)$. The next section explains the Hilbert space Gaussian process (HSGP) approximation to build an effective proxy of $\imse(\bft)$.
\begin{algorithm}
    \caption{Sequential sampling with $\imse(\bft)$} \label{alg:sequential_with_exact_imse}
    \renewcommand{\algorithmicrequire}{\textbf{Input:}}
    \begin{algorithmic}
        \Require Initial data $\{(\bfx_n,y_n)\}_{n=1}^N$ and an iteration budget $\texttt{step}$
        \State Fit a GP model using \eqref{eq:log-likelihood}--\eqref{eq:mle_numerical}
        \For{$i=1,\ldots,\texttt{step}$}
            \State Compute $\imse(\bft)$ for candidate $\bft\in\Omega$ (typically via numerical integration)
            \State $\bfx_{N+i}\gets \arg\max_{\bft\in\Omega}\imse(\bft)$
            \State Observe $y_{N+i}=f(\bfx_{N+i}) + \epsilon_{N+i}$
            \State Update the GP fit using $\{(\bfx_n,y_n)\}_{n=1}^{N+i}$
        \EndFor
        \State \Return fitted GP model based on $\{(\bfx_n,y_n)\}_{n=1}^{N+\texttt{step}}$
    \end{algorithmic}
\end{algorithm}

\subsection{HSGP-based approximation for IMSE} \label{subsec:HSGP-IMSE}
The method of $\hsgp$ is first proposed by \citet{solin:2020} as a computationally efficient reduced-rank GP regression method for stationary kernels, and its empirical performance has been investigated in detail in \citet{riutort:2023}. In this section, we first introduce the main technical steps behind the idea of $\hsgp$ approximation of stationary kernels, and then explain how to apply it to sequential designs using an induced approximation of $\imse(\bft)$.

Let $m\in\NN$ and define $[m]^d=\{1,\ldots,m\}^d$. Fix $L>0$ and define the sine basis on $(-L,L)^d$ by
\begin{align} \label{eq:hsgp_sin_basis}
    \phi_{\bfj}(\bfx)
    := L^{-d/2}\prod_{k=1}^d \sin \left(\frac{\pi j_k}{2L}(x_k+L)\right),
    \qquad \bfj=[j_1,\ldots,j_d]^\T\in\ZZ_+^d.
\end{align}
These functions are the eigenfunctions of the Laplace operator on $(-L,L)^d$ with Dirichlet boundary conditions:
\begin{align} \label{eq:laplace_pde}
    \begin{cases}
    -\Delta \phi_{\bfj}(\bfx)=\lambda_{\bfj}\phi_{\bfj}(\bfx), & \bfx\in(-L,L)^d,\\
    \phi_{\bfj}(\bfx)=0, & \bfx\in\partial(-L,L)^d.
    \end{cases}
\end{align}
The corresponding eigenvalues are given by
\begin{align}
    \lambda_{\bfj}
    =
    \left(\dfrac{\pi}{2L}\right)^2 \|\bfj\|_2^2, \qquad \bfj=[j_1,\ldots,j_d]^\T\in\ZZ_+^d.
    \nonumber
\end{align}
The spectral density $S(\omega)$ of shift-invariant kernel $k$ is defined as the Fourier transform:
\begin{align} \label{eq:spec_def}
    S(\omega) = \int_{\RR^d} k(\bfh)  \exp \left(- i \omega^\T \bfh \right)  \ud \bfh. 
\end{align}
From formula (4) to (20) in \citet{solin:2020}, for a shift-invariant kernel with spectral density $S(\omega)$, defining the operator $T$ as follows:
\begin{align} \label{eq:kernel_operator_def}
    T f = \int_{\RR^d} k(\cdot, \bfx) f(\bfx) \ud \PP(\bfx),\qquad \forall f \in L^2(\Omega).
\end{align}
For a heuristic argument, \citet{solin:2020} assumed that $S(\omega) = S(\sqrt{\|\omega\|_2^2})$ is an analytic function of $\|\omega\|_2^2$ with the Taylor series expansion:
\begin{align} 
    S(\omega) = \sum_{k=0}^\infty a_k \|\omega\|_2^{2k}. \nonumber
\end{align}
This is true for infinitely smooth kernels such as the Gaussian kernel. Since $\mathcal{F}(-\Delta f) (\omega) = \|\omega\|_2^2 \mathcal{F}(f) (\omega)$, where $\mathcal{F}(f)$ is the Fourier transform of a function $f$, it follows that
\begin{align} \label{eq:kernel_operator_identity}
    T = a_0 I + \sum_{k=1}^\infty a_k (-\Delta)^k.
\end{align}
Then based on the eigen-decomposition specified by \eqref{eq:hsgp_sin_basis} and \eqref{eq:laplace_pde}, plugging \eqref{eq:kernel_operator_identity} into \eqref{eq:kernel_operator_def} yields the HSGP approximation proposed by \citet{solin:2020}, 
\begin{align*}
    k(\bfx, \bfx^\prime) 
    \approx 
    \sum_{\bfj} S\left(\sqrt{\lambda_\bfj}\right) \phi_\bfj(\bfx) \phi_\bfj(\bfx^\prime),
\end{align*}
where $\lambda_\bfj$, $\phi_\bfj$ are defined as the solutions to the eigenvalue problem of PDE \eqref{eq:laplace_pde}.
Thus, the infinite-series HSGP representation is
\begin{align} \label{eq:kernel_hsgp_infinite}
    \hat k_\infty(\bfx,\bfx')
    := \sum_{\bfj\in\ZZ_+^d}
    S \left(\frac{\pi\bfj}{2L}\right) 
    \phi_{\bfj}(\bfx) \phi_{\bfj}(\bfx').
\end{align}
In practice, truncating to $[m]^d$ yields the finite approximation
\begin{align} \label{eq:kernel_hsgp_finite}
    \hat k_m(\bfx,\bfx')
    := \sum_{\bfj\in[m]^d}
    S \left(\frac{\pi\bfj}{2L}\right) 
    \phi_{\bfj}(\bfx) \phi_{\bfj}(\bfx').
\end{align}
For later use, we provide the spectral densities of the Gaussian and Mat\'ern kernels under the Fourier convention stated in \eqref{eq:spec_def}:
\begin{itemize}
    \item \textbf{Gaussian kernel spectral density:}
    \begin{align} \label{eq:spec_gaussian}
        S(\omega)=\sigma^2 (2\pi)^{d/2}\ell^d \exp \left(-\frac{\ell^2\|\omega\|_2^2}{2}\right).
    \end{align}
    \item \textbf{Mat\'ern kernel spectral density:}
    \begin{align} \label{eq:spec_matern}
        S(\omega)
        = \sigma^2 \frac{2^d \pi^{d/2}\Gamma(\nu+d/2)}{\Gamma(\nu)}
        \left(\frac{2\nu}{\ell^2}\right)^{\nu}
        \left(\frac{2\nu}{\ell^2}+\|\omega\|_2^2\right)^{-(\nu+d/2)}.
    \end{align}
\end{itemize}

\paragraph{HSGP approximation for IMSE.} Based on this HSGP heuristic, replacing $k$ in \eqref{eq:imse_int} by $\hat k_m$ leads to the HSGP approximation of $\imse(\bft)$:
\begin{align} \label{eq:imse_hsgp_finite}
    \widehat{\imse}_m(\bft)
    := \frac{1}{P_{N, \eta}^2(\bft) + \eta}
    \int_{\Omega}\left[\hat k_m(\bfx,\bft)-\hat\bfk_{N,m}^\T(\bfx) \left(\K_N + \eta \I_N \right)^{-1} \bfk_N(\bft)\right]^2  \ud \PP(\bfx),
\end{align}
where $\hat\bfk_{N,m}(\bfx):=[\hat k_m(\bfx,\bfx_1),\ldots,\hat k_m(\bfx,\bfx_N)]^\T$. We emphasize that we replace only the kernel functions that involve the integrated variable $\bfx$ with the HSGP approximation.

One can immediately see the major advantage of HSGP approximation when evaluating integrals that involve kernel functions. Since the HSGP approximation $\hat k_m(\bfx,\bfx')$ in \eqref{eq:kernel_hsgp_finite} depends only on trigonometric functions, the integral of $\hat k_m(\bfx,\bfx')$ on rectangle domains such as $\Omega=(-B,B)^d$ has a closed form, because the trigonometric functions can be integrated in closed form.  Therefore, as long as the spectral densities of the kernel $k(\bfx,\bfx')$ can be evaluated in closed form, the HSGP approximation enables efficient computation of kernel integrals. Moreover, the two hyperparameters $L$ and $m$ can be chosen to ensure that the approximation error is small. In general, larger values of $L$ lead to smaller errors.

\begin{remark}[Choice of $L$ and boundary effects]
In practice, the HSGP basis $\{\phi_{\bfj}\}$ is orthogonal normal basis on a larger domain $(-L,L)^d$ with $L>B$ and $\Omega = (-B, B)^d \subset (-L, L)^d$, although the design and integration domain is still $\Omega=(-B,B)^d$. This padding pushes the Dirichlet boundary $\partial(-L,L)^d$ away from $\Omega$, mitigating boundary artifacts of the truncated series \eqref{eq:kernel_hsgp_finite} inside $\Omega$.
\end{remark}

\begin{theorem} \label{thm:HSGP-IMSE_acq}
Suppose $\PP$ is the Lebesgue measure and $\Omega=(-B,B)^d$. Suppose $\phi_\bfj$ is the HSGP basis defined in \eqref{eq:hsgp_sin_basis}. Let $M=m^d$, and index the truncated basis $\{\phi_{\bfj}\}_{\bfj\in[m]^d}$ in any fixed order so that
$
\phi(\bfx) := \bigl[\phi_1(\bfx),\ldots,\phi_M(\bfx)\bigr]^\T.
$
Define the design matrix $\Phi = \left[\phi_i(\bfx_n) \right] \in\RR^{N\times M}, i \in \{1, \ldots, M\}, n \in \{1, \ldots, N\}$. Define the spectral matrix and the Gram matrix
\begin{align}
    \W := \diag \left[S \left(\frac{\pi\bfj}{2L}\right)\right]_{\bfj\in[m]^d}, \quad 
    \G_d := \int_\Omega \phi(\bfx)\phi^\T(\bfx) \ud \PP(\bfx). \nonumber
\end{align}
Let
\begin{align}
    \bfh(\bft) 
    =
    \phi(\bft)-\Phi^\T \left(\K_N + \eta \I_N \right)^{-1}\bfk_N(\bft). \nonumber
\end{align}
Then $\widehat{\imse}_m(\bft)$ admits the closed form
\begin{align} \label{eq:HSGP-IMSE_acq}
    \widehat{\imse}_m(\bft)
    = \frac{\bfh^\T(\bft) \W \G_d \W \bfh(\bft)}{P_{N, \eta}^2(\bft) + \eta} = \frac{\bfh^\T(\bft) \W \left[\G_1^{\otimes d} \right]\W \bfh(\bft)}{P_{N, \eta}^2(\bft) + \eta},
\end{align}
where 
\begin{align}
    [\G_1]_{p,q}
    &=
    \int_{-B}^{B}\phi_p(x)\phi_q(x) \ud x \nonumber \\
    &= \begin{cases}
       \frac{B}{L} - \frac{1}{2\pi p} \left[\sin(\frac{\pi p}{L}(L+B)) - \sin(\frac{\pi p}{L}(L-B)) \right], & p = q, \\
       \frac{1}{\pi} \left[\frac{1}{p-q} \sin(\frac{\pi}{2 L}(p-q)s) - \frac{1}{p+q} \sin(\frac{\pi}{2 L}(p+q)s) \right] \Big |_{s=L-B}^{s=L+B}, & p \neq q,
    \end{cases} \nonumber.
\end{align}
\end{theorem}

\paragraph{Complexity analysis.} Evaluating $\widehat{\imse}_m(\bft)$ involves two main steps:  
(i) linear algebra for calculating $\left(\K_N + \eta \I_N \right)^{-1}\bfk_N(\bft)$, typically obtained via a Cholesky factorization, at the cost $O(N^3)$; and  
(ii) multiplication by $\G_d$ and the diagonal weight matrix $\W$. 
The Kronecker structure $G_d = \G_1^{\otimes d}$ allows multiplication by $\G_d$ via mode-wise products at cost $O(d m M)$, rather than standard matrix multiplication at cost $O(m^{2d})$. This Kronecker structure is especially beneficial when $d>1$, reducing the cost relative to forming $\G_d$ explicitly.  Hence step (ii) requires a cost of
$
O\bigl((N+dm)M\bigr).
$
In Section~\ref{subsec:HSGP.SD.rate}, we will show that fast convergence of $\|\widehat{\imse}_m(\bft)-\imse(\bft)\|_{L^q}, q \in [1, \infty)$ requires
\begin{align*}
M \asymp N^{\frac{\nu-d/(2q)}{\nu}}(\log N)^d =o(N).
\end{align*}
Consequently, the computational cost of step (ii) is 
$
O(NM)
=
o(N^2).
$
By contrast, \citet{binois:2019} rewrite the numerator of $\imse(\bft)$ as
\begin{align*}
    w(\bft, \bft) - 2 \bfk_N^\T(\bft) \left(\K_N + \eta \I_N \right)^{-1} \bfw_N(\bft) + \bfk_N^\T(\bft) \left(\K_N + \eta \I_N \right)^{-1} \W_N \left(\K_N + \eta \I_N \right)^{-1} \bfk_N(\bft),
\end{align*}
where $w(\bfx, \bfx^\prime) = \int_{\Omega} k(\bfx, \bfz) k(\bfz, \bfx^\prime) \ud \PP(\bfz)$, $\bfw_N(\bft) = \left[w(\bft, \bfx_1), \ldots, w(\bft, \bfx_N) \right]^\T$, and $\W_N = [W_{ij}]_{i, j = 1}^N$. They pre-compute each entry $W_{ij}$ as a function of $\bfx_i$ and $\bfx_j$. This approach also consists of two steps:  
(i) solving for $\left(\K_N + \eta \I_N \right)^{-1}\bfk_N(\bft)$ and $\left(\K_N + \eta \I_N \right)^{-1} \bfw_N(\bft)$;  
(ii) performing the linear algebra required for $w(\bft, \bft) - 2 \bfk_N^\T(\bft) \K_N^{-1} \bfw_N(\bft) + \bfk_N^\T(\bft) \K_N^{-1} \W_N \K_N^{-1} \bfk_N(\bft)$, which can be done via Cholesky decomposition. Step (i) again has the cost $O(N^3)$. Step (ii) requires $O(N^2+N)=O(N^2)$ operations to evaluate matrix products. From this perspective, our HSGP-IMSE method is computationally more efficient because it requires only $o(N^2)$ cost.

The above discussion assumes that the kernel parameters are re-estimated at each iteration. In that case, the kernel must be re-evaluated at every iteration for all design points $\{\bfx_n\}_{n=1}^N$, and the $O(N^3)$ cost of step (i) cannot in general be avoided. If the parameters are instead held fixed during an iteration, then step (i) can be updated recursively via the block matrix inverse formula: $\K_{N+1}^{-1}$ can be obtained from $\K_N^{-1}$ in $O(N^2)$ time. Under this setting, the HSGP-IMSE method remains faster.

\subsection{$\gamma$-stabilizing strategy} \label{subsec:gamma_greedy}
The maximization of $\widehat{\imse}_m(\bft)$ may become numerically unstable when candidates approach existing design points, as is illustrated in Example 5 in \citet{wenzel:2021}. A stabilizing strategy is to impose a separation constraint during acquisition maximization:
\begin{align} 
&\arg\max_{\bft\in\Omega_\gamma}\widehat{\imse}_m(\bft), \nonumber \\
\text{where } \quad  & \Omega_\gamma
    := \Omega \cap \left\{\bft:\dist(\bft,\X_N)\ge \gamma h_N\right\}, \nonumber \\
    & \dist(\bft,\X_N)
    := \min_{\bfx_k\in\X_N}\|\bft-\bfx_k\|_2, \nonumber \\
    & h_N    := \sup_{\bfx\in\Omega}\min_{\bfx_k\in\X_N}\|\bfx-\bfx_k\|_2, \nonumber 
\end{align}
and $\gamma\in(0,1)$ is a fixed hyper-parameter. Here $h_N$ is the fill distance of $\X_N$ in $\Omega$. The optimization over the space $\Omega_{\gamma}$ prevents the next point $\bfx_{N+1}$ from being too close to any current points in $\X_N$. Thus, the constraint encourages exploration by preventing overly local refinement. This $\gamma$-stabilizing strategy also helps to ensure asymptotical quasi-uniform property, which is shown in the following Theorem~\ref{thm:quasi_uniform_points}. 
\begin{theorem}[Quasi-uniformity under $\gamma$-stabilizing sampling] \label{thm:quasi_uniform_points}
Define
$
q_N := \frac{1}{2}\min_{i\neq j}\|\bfx_i-\bfx_j\|_2,
$
which is the separation distance of $\X_N=\{\bfx_1,\ldots,\bfx_N\} \subset \Omega$. Fix $\gamma\in(0,1)$ and suppose $\bfx_{N+1}$ is selected from $\Omega_\gamma$. Then for any $N\ge 2$,
\begin{align} 
    \frac{h_N}{q_N} \le 2\gamma^{-1}. \nonumber
\end{align}
\end{theorem}

Theorem~\ref{thm:quasi_uniform_points} shows that the points $\bfx_1,\bfx_2,\ldots$ obtained from the $\gamma$-stabilizing strategy are guaranteed to be quasi-uniform in $\Omega$. 
Based on the stabilized criterion $\widehat{\imse}_m(\bft)$, the resulting sequential procedure is summarized in Algorithm~\ref{alg:sequential_with_HSGP-IMSE}.
\begin{algorithm}
    \caption{Sequential sampling with $\widehat{\imse}_m(\bft)$} \label{alg:sequential_with_HSGP-IMSE}
    \renewcommand{\algorithmicrequire}{\textbf{Input:}}
    \begin{algorithmic}
        \Require Initial data $\{(\bfx_n,y_n)\}_{n=1}^N$ and an iteration budget $\texttt{step}$
        \State Fit the GP model using \eqref{eq:log-likelihood}--\eqref{eq:mle_numerical}
        \State Select initial HSGP parameters $(m,L)$
        \For{$i=1,\ldots,\texttt{step}$}
            \State Compute $\widehat{\imse}_m(\bft)$ based on the current data  $\{(\bfx_n,y_n)\}_{n=1}^{N+i-1}$ using Theorem~\ref{thm:HSGP-IMSE_acq} 
            \State $\bfx_{N+i}\gets \arg\max_{\bft\in\Omega_\gamma}\widehat{\imse}_m(\bft)$
            \State Observe $y_{N+i}=f(\bfx_{N+i})+ \epsilon_{N+i}$
            \State Update the GP fit using $\{(\bfx_n,y_n)\}_{n=1}^{N+i}$
            \State Optionally update $(m,L)$
        \EndFor
        \State \Return fitted GP model based on $\{(\bfx_n,y_n)\}_{n=1}^{N+\texttt{step}}$
    \end{algorithmic}
\end{algorithm}

\vspace{5mm}

\section{Theory} \label{sec:theory}
This section establishes approximation guarantees for the $\hsgp$-based surrogate $\widehat{\imse}_m(\bft)$ to the exact acquisition function $\imse(\bft)$. The analysis is organized in two steps.
Section~\ref{subsec:HSGP.rate} studies the truncation error induced by the finite-dimensional $\hsgp$ representation of the kernel. Explicit constants are provided in the resulting convergence bounds, so that the dependence on the truncation level $m$, the domain parameter $L, B$, and the kernel parameters $\sigma^2, \ell$ is transparent.
Building on these kernel-level estimates, Section~\ref{subsec:HSGP.SD.rate} analyzes the induced error in the acquisition function for nugget/noise-free cases, i.e., $\eta=0$ in \eqref{eq:model_gpr}. In particular, an upper bound is derived for the approximation error of the numerator in $\widehat{\imse}_m(\bft)$, together with a uniform lower bound for the exact denominator (the squared power function). Combining these two ingredients yields a convergence rate for the overall approximation error
$
    \widehat{\imse}_m(\bft)-\imse(\bft),
$
thereby quantifying the accuracy of the $\hsgp$-based approximation across the design domain.
Section~\ref{subsec:nugget} further derives the upper bound in the presence of nugget, i.e., $\eta>0$ in \eqref{eq:model_gpr}. Our results include both the cases where $\eta$ is allowed to depend on $N$ and where $\eta$ is of the constant order.

\subsection{Sharp kernel approximation rates for Hilbert space Gaussian process} \label{subsec:HSGP.rate}

This subsection quantifies the accuracy of the Hilbert space Gaussian process ($\hsgp$) approximation for shift-invariant kernels on a bounded design domain. The approximation error is decomposed into two conceptually distinct components: an \emph{aliasing error} induced by restricting the infinite-domain kernel to the finite computational domain $(-L,L)^d$, and a \emph{truncation error} induced by retaining only finitely many basis functions in the HSGP expansion. While \citet{solin:2020} provides error bounds that decay polynomially with respect to $L$, we provide improved bounds that decay exponentially in $L$ with explicit constants. 

\paragraph{Error decomposition.} Recall the infinite and truncated HSGP approximations $\hat k_\infty$ and $\hat k_m$ defined in \eqref{eq:kernel_hsgp_infinite}--\eqref{eq:kernel_hsgp_finite}. For any $\bfx,\bft\in\Omega$, define the kernel approximation error as
\begin{align} \label{eq:error_def_kernel_hsgp}
    \hat k_m(\bfx, \bft) - k(\bfx, \bft)
    &=
    \underbrace{\left[\hat k_\infty(\bfx, \bft) - k(\bfx, \bft)\right]}_{\text{aliasing error}}
    -
    \underbrace{\left[\hat k_\infty(\bfx, \bft) - \hat k_m(\bfx, \bft)\right]}_{\text{truncation error}} .
\end{align}
The aliasing term reflects the discrepancy between the target kernel on $\Omega$ and the infinite HSGP series constructed from eigenfunctions on $(-L,L)^d$. The truncation term measures the tail of the infinite series beyond the finite index set. In both cases, uniform (supremum) bounds over $\Omega\times\Omega$ are derived, as these are convenient for subsequent error propagation to acquisition functions.

\paragraph{Assumption 1.} $\Omega = (-B, B)^d \subset \RR^d, B > 0$, which is a bounded domain satisfying Lipschitz boundary and interior cone conditions. $\X_N := \{\bfx_n\}_{n=1}^N \subset \Omega$ consists of distinct points.

This regularity condition is standard in approximation theory and ensures that constants appearing in uniform bounds can be controlled in terms of geometric properties of $\Omega$. The computational domain parameter $L$ is chosen such that $L>B$, so that $\Omega$ is strictly contained in the expansion domain $(-L,L)^d$.

\begin{theorem}[Aliasing and truncation errors for the Gaussian kernel] \label{thm:gaussian_kernel_hsgp_error}
Consider the Gaussian kernel defined in \eqref{eq:kernel_gaussian} with spectral density given in \eqref{eq:spec_gaussian}. Suppose Assumption 1 holds and $L>B$. Then the aliasing error admits the uniform bound
\begin{align} 
    \sup_{\bfx, \bft \in \Omega} \bigl|\hat k_\infty(\bfx, \bft) - k(\bfx, \bft)\bigr|
    \leq 
    \sigma^2 (2^d + 2d-1)\left(3 + \frac{\sqrt{2\pi} \ell}{4L}\right)^{d} \exp \left(-\frac{2(L-B)^2}{\ell^2} \right). \nonumber
\end{align}
Moreover, if $L > \ell \sqrt{\pi/2}$, the truncation error admits the uniform bound
\begin{align} 
    \sup_{\bfx, \bft \in \Omega} \bigl|\hat k_\infty(\bfx, \bft) - \hat k_m(\bfx, \bft)\bigr|
    \leq 
    \frac{4\sqrt{2}}{\pi^{3/2}\ell}  3^{d-1} d \sigma^2  \frac{L}{m} 
    \exp \left(-\frac{\pi^2 \ell^2 m^2}{8L^2}\right). \nonumber
\end{align}
\end{theorem}

This theorem shows that the aliasing error for the Gaussian kernel decays exponentially fast as $L$ increases, consistent with the rapid tail decay of the Gaussian kernel. For a fixed $L$, the truncation error also decays exponentially fast in $m$, reflecting the super-polynomial smoothness of the Gaussian kernel. 

\begin{corollary}
    Under the assumptions in Theorem~\ref{thm:gaussian_kernel_hsgp_error}, equating the aliasing and truncation terms by setting $m =O(L^2)$, the kernel approximation error \eqref{eq:error_def_kernel_hsgp} admits the convergence rate
    \begin{align}
        \sup_{\bfx, \bft \in \Omega} \bigl|\hat k_\infty(\bfx, \bft) - k(\bfx, \bft)\bigr| \lesssim \exp \left( -\frac{\pi}{2} m \right). \nonumber
    \end{align}
\end{corollary}

\medskip

\begin{theorem}[Aliasing and truncation errors for the Mat\'ern kernel] \label{thm:matern_kernel_hsgp_error}
Consider the \\ Mat\'ern kernel defined in \eqref{eq:kernel_matern} with spectral density given in \eqref{eq:spec_matern}. Suppose Assumption 1 holds and
\begin{align*}
    L > \max\left\{\frac{1}{2}(B + \sqrt{2 d \nu} \ell), ~\frac{\sqrt{d} \ell }{2\sqrt{2\nu}}\log 3\right\}.
\end{align*}
Then the aliasing error admits the uniform bound
\begin{align} 
    \sup_{\bfx, \bft \in \Omega} \bigl|\hat k_\infty(\bfx, \bft) - k(\bfx, \bft)\bigr|
    \leq 
    \sigma^2 (d + 2^d - 1) \frac{2^{d+\nu+2}}{\Gamma(\nu)} \nu^\nu K_\nu(4\nu) \exp\left(2\nu +\frac{\sqrt{2\nu}}{\sqrt{d} \ell}B -\frac{\sqrt{2\nu}}{\sqrt{d} \ell} L\right). \nonumber
\end{align}
Moreover, the truncation error admits the uniform bound
\begin{align} 
    \sup_{\bfx, \bft \in \Omega} \bigl|\hat k_\infty(\bfx, \bft) - \hat k_m(\bfx, \bft)\bigr|
    \leq  
    \sigma^2 2^{2\nu+d+1} d \frac{\Gamma(\nu + d/2)}{\Gamma(\nu)} \pi^{-(2\nu+d/2)} \left(\frac{2\nu}{\ell^2}\right)^{\nu} \beta(d, \nu) 
    (L/m)^{2\nu}, \nonumber
\end{align}
with 
\begin{align}
    \beta(d, \nu) = \begin{cases}
        \frac{1}{2\nu}, & \text{for } d = 1, \\
        \left(4+\frac{2}{2\nu+d-1}\right) \beta(d-1, \nu), & \text{for } d > 1.
    \end{cases} \nonumber
\end{align}
\end{theorem}

The Mat\'ern bounds exhibit a qualitative difference from the Gaussian case. While the aliasing error still decays exponentially in $L$, the truncation error scales polynomially in $m$ at rate $(L/m)^{2\nu}$. This dependence reflects the finite smoothness of Mat\'ern sample paths, controlled by the parameter $\nu$. 

\begin{corollary}
    Under the assumptions in Theorem~\ref{thm:matern_kernel_hsgp_error}, by setting $L \asymp \sqrt{2\nu d} \ell \log m$, the kernel approximation error \eqref{eq:error_def_kernel_hsgp} admits the convergence rate
    \begin{align}
        \sup_{\bfx, \bft \in \Omega} \bigl|\hat k_\infty(\bfx, \bft) - k(\bfx, \bft)\bigr| \lesssim \left(\frac{\log m}{m}\right)^{2\nu}. \nonumber
    \end{align}
\end{corollary}

\begin{remark}
    Theorems \ref{thm:gaussian_kernel_hsgp_error} and \ref{thm:matern_kernel_hsgp_error} fill the gap in the literature of HSGP approximation convergence rate. While \citet{solin:2020} provide a convergence analysis and conclude that $|\hat k_m(\bfx, \bft) - k(\bfx, \bft)| \leq O(dL^{-1} + m^{-2\nu/d})$, this error bound is not sharp, and the associated constants are not given in explicit form.
\end{remark}

\paragraph{Implications for selecting $(L,m)$.} Theorems~\ref{thm:gaussian_kernel_hsgp_error} and~\ref{thm:matern_kernel_hsgp_error} separate the roles of the two tuning parameters. Increasing $L$ reduces aliasing error, whereas increasing $m$ reduces truncation error. The explicit forms of the bounds will be used in Section~\ref{subsec:HSGP.SD.rate} to control the numerator error in $\widehat{\imse}_m(\bft)$ and to derive convergence rates for the approximation of $\imse(\bft)$.

\subsection{Error rates for HSGP-based IMSE approximation without noise/nugget} \label{subsec:HSGP.SD.rate}

This subsection derives quantitative error bounds for the HSGP-based acquisition $\widehat{\imse}_m(\bft)$ when the underlying GP prior uses a Mat\'ern kernel. The noise/nugget term $\eta$ is assumed to be 0. Let $P_N(\bft) = P_{N, 0}(\bft)$, which is the standard power function. For any $\bft \in \Omega_\gamma$,
\begin{align*}
    \imse(\bft) = \frac{1}{P_N^2(\bft)}\int_{\Omega}
    \Bigl[k(\bfx,\bft) - \bfk_N^\T(\bfx) \K_N^{-1} \bfk_N(\bft)\Bigr]^2
     \ud \PP(\bfx).
\end{align*}
Let $\|g(\bft) \|_{L^q} = \|g(\bft) \|_{L^q(\Omega_\gamma, \PP, \bft)} = 
\left[\int_{\Omega_\gamma} [g(\bft)]^q \ud \PP(\bft) \right]^{1/q}$ denote the $L^q$ norm with respect to $\bft$, the goal is to control the discrepancy $\|\widehat{\imse}_m(\bft)-\imse(\bft)\|_{L^q}$ for $q \in [1, \infty]$
over candidate locations $\bft\in\Omega_\gamma$, and to translate kernel-level HSGP approximation rates into acquisition-level rates relevant for sequential design. For instance, the following two metrics,
\begin{align*}
    \|\widehat{\imse}_m(\bft)-\imse(\bft)\|_{L^1}, \quad \|\widehat{\imse}_m(\bft)-\imse(\bft)\|_{L^\infty}.
\end{align*}
represent the expected error and the worst case error. 

\paragraph{Error representation.} Starting from Theorem~\ref{thm:imse_acq} and its HSGP analogue, the pointwise error can be written as
\begin{align}\label{eq:error_def_imse_hsgp}
    &~\quad \widehat{\imse}_m(\bft) - \imse(\bft) \nonumber \\
    &= \frac{1}{P_N^2(\bft)}
    \Biggl\{
    \int_{\Omega}
    \Bigl[\hat{k}_m(\bfx,\bft) - \bfk_N^\T(\bft) \K_N^{-1} \hat{\bfk}_{N,m}(\bfx)\Bigr]^2
     \ud \PP(\bfx)
    -\int_{\Omega}
    \Bigl[k(\bfx,\bft) - \bfk_N^\T(\bfx) \K_N^{-1} \bfk_N(\bft)\Bigr]^2
     \ud \PP(\bfx)
    \Biggr\} \nonumber\\
    &= \frac{1}{P_N^2(\bft)} \int_\Omega \left\{[\hat{k}_m(\bfx,\cdot) - \Pi_{\X_N} \hat{k}_m(\bfx,\cdot)]^2 - [k(\bfx,\cdot) - \Pi_{\X_N} k(\bfx,\cdot)]^2 \right\} \ud \PP(\bfx),
\end{align}
where $[\Pi_{\X_N} f](\cdot) := \bfk_N^\T(\cdot) \K_N^{-1} \left[f(\bfx_1), \ldots, f(\bfx_N)\right]^\T$ is the projection of $f$ onto the kernel space spanned by $\bfk(\bfx_1,\cdot),\ldots,\bfk(\bfx_N,\cdot)$.

The structure in \eqref{eq:error_def_imse_hsgp} makes clear that the approximation error is driven by two ingredients:
(i) the error in approximating the \emph{numerator} integral (a squared residual integrated over $\Omega$), and
(ii) the instability of the \emph{denominator} $P_N^2(\bft)$, which can become small when $\bft$ approaches the current design $\X_N$.
The analysis therefore proceeds by bounding the numerator error from above and bounding the denominator $P_N^2(\bft)$ from below on the constrained set $\Omega_{\gamma}$ as specified in the $\gamma$-stabilizing strategy in Section~\ref{subsec:gamma_greedy}. We first summarize the necessary assumptions.

\paragraph{Assumption 2.} $k$ is a Mat\'ern kernel defined in \eqref{eq:kernel_matern} with fixed smoothness $\nu > d/2$. Let $\tau = \nu + d/2$.
\paragraph{Assumption 3.} Let $h_N$ denote the fill distance of $\{\bfx_n\}_{n=1}^N$. There exists a constant $h_0>0$ (determined by the interior cone condition) such that $h_N\leq h_0$.

\vspace{2mm}

Assumption 2 specifies the condition on the smoothness of kernel. It ensures that the Mat\'ern RKHS is continuously embedded into a suitable Sobolev space, so that the kernel sections are sufficiently smooth and the interpolation error admits the standard fill-distance bounds used later.
Assumption 3 means that the design points are dense enough in $\Omega$. 
The threshold $h_0$ depends only on the geometry of the domain, and condition $h_N \le h_0$ guarantees standard scattered-data approximation and kernel interpolation estimations.

\begin{theorem} \label{thm:error_imse_hsgp_numerator}
Suppose that Assumptions 1-3 hold. Suppose
\begin{align}
    m > L > 
    B + \frac{1}{2} \max\left\{1, ~2\sqrt{2\nu} \ell, ~\frac{4(\tau+1) \ell}{\sqrt{2\nu}}\right\}. \nonumber
\end{align}
Then 
\begin{align}
    &~\quad \left\| \int_\Omega \left\{[\hat{k}_m(\bfx,\cdot) - \Pi_{\X_N} \hat{k}_m(\bfx,\cdot)]^2 - [k(\bfx,\cdot) - \Pi_{\X_N} k(\bfx,\cdot)]^2 \right\} \ud \PP(\bfx) \right\|_{L^q} \nonumber\\
    &\lesssim h_N^{2\nu + \frac{d}{q}} \left[ \exp \left(-\frac{\sqrt{2\nu}}{2\sqrt{d} \ell} L\right) + (L/m)^{\nu} \right]^2 
    + h_N^{3\nu + \frac{d}{2q}} \left[\exp \left(-\frac{\sqrt{2\nu}}{2\sqrt{d} \ell} L\right) + (L/m)^{\nu} \right]. \nonumber
\end{align}
\end{theorem}

Theorem~\ref{thm:error_imse_hsgp_numerator} shows that the numerator error inherits two distinct HSGP kernel approximation mechanisms: an \emph{aliasing component} that decays exponentially in $L$ and a \emph{truncation component} that decays polynomially in $m$ at rate $(L/m)^\nu$ for Mat\'ern kernels.

\medskip

\begin{theorem} \label{thm:lower_bound_power_function}
Suppose that Assumptions 1-3 hold. With the $\gamma$-stabilizing strategy described in Section~\ref{subsec:gamma_greedy}, for any $\bft \in \Omega_{\gamma}$ with $\dist(\bft,\X_N)\ge \gamma h_N$,
\begin{align}
    P_N^2(\bft) \gtrsim h_N^{2\nu}. \nonumber
\end{align}
\end{theorem}

Theorem~\ref{thm:lower_bound_power_function} provides the key stability ingredient: once the candidate $\bft$ is kept a fixed fraction of the fill distance away from the current design, the denominator $P_N^2(\bft)$ cannot collapse to zero. This is precisely the role of the $\gamma$-stabilizing restriction introduced in Section~\ref{subsec:gamma_greedy}. 

\medskip

\begin{corollary}[Upper bound for $\widehat{\imse}_m(\bft)$ in noise/nugget free case] \label{cor:global_bound_imse_hsgp}
Suppose conditions in Theorems~\ref{thm:error_imse_hsgp_numerator} and~\ref{thm:lower_bound_power_function} hold. With the same notations,
\begin{align}
    \bigl\|\widehat{\imse}_m(\bft) - \imse(\bft)\bigr\|_{L^q}
    &\lesssim
    h_N^{\frac{d}{q}} \left[ \exp \left(-\frac{\sqrt{2\nu}}{2\sqrt{d} \ell} L\right) + (L/m)^{\nu} \right]^2  \nonumber\\
    &\quad + h_N^{\nu + \frac{d}{2q}} \left[\exp \left(-\frac{\sqrt{2\nu}}{2\sqrt{d} \ell} L\right) + (L/m)^{\nu} \right]. \nonumber
\end{align}
\end{corollary}

Corollary~\ref{cor:global_bound_imse_hsgp} is a direct result of combining Theorem~\ref{thm:error_imse_hsgp_numerator} and \ref{thm:lower_bound_power_function}. It combines (i) control of the numerator perturbation and (ii) a uniform lower bound on the power function over $\gamma$-stabilizing feasible candidates, yielding an explicit acquisition-level convergence rate. In particular, the bound separates the dependence on the design resolution ($h_N$) from the dependence on the HSGP truncation parameters $(L,m)$, thereby enabling principled tuning.

\begin{remark}[Guidance for choosing $(L,m)$] \label{rem:choice_L_m}
Corollary~\ref{cor:global_bound_imse_hsgp} suggests balancing the aliasing and truncation errors by enforcing 
$
(L/m)^\nu \approx \exp \left(-\frac{\sqrt{2\nu}}
    {2\sqrt{d} \ell}L\right).
$
Let $L = \sqrt{2d\nu}\ell \log m$, then 
$
(L/m)^\nu 
\asymp \exp \left(-\frac{\sqrt{2\nu}}{2\sqrt{d} \ell}L\right) 
\asymp \big(\frac{\log m}{m}\big)^\nu.
$
Under the quasi-uniform design in Theorem~\ref{thm:quasi_uniform_points}, we have $h_N = O(N^{-1/d})$. Consequently, the error bound turns out to be
\begin{align*}
    N^{-\frac{1}{q}} \left(\frac{\log m}{m}\right)^{2\nu} 
    +
    N^{-\frac{1}{2q} - \frac{\nu}{d}} \left(\frac{\log m}{m}\right)^\nu.
\end{align*}
In principle, we require that $m \lesssim N^{1/d}$ such that $m^d \lesssim N$. Therefore, by taking $m = O(N^{1/d})$, the optimal error bound is as follows
\begin{align*}
    N^{-\frac{1}{q} - \frac{2\nu}{d}} (\log N)^{2\nu} 
    +
    N^{-\frac{1}{2q} - \frac{2\nu}{d}} (\log N)^\nu
    =
    \begin{cases}
        N^{-\frac{1}{2q} - \frac{2\nu}{d}} (\log N)^\nu, & q < \infty, \\
        N^{- \frac{2\nu}{d}} (\log N)^{2\nu}, & q = \infty.
    \end{cases}.
\end{align*}
\end{remark}

\begin{remark}
    By the Cauchy-Schwarz inequality, $\imse(\bft) \leq \int_\Omega P_N^2(\bfx) \ud \PP(\bfx) \lesssim N^{-2\nu/d}$, since $\sup_{\bfx\in \Omega} P_N^2(\bfx) \lesssim h_N^{2\nu} = O(N^{-2\nu/d})$ under the quasi-uniform design by Theorem 11.11 in \citet{wendland:2005}. By triangle inequality, if we also have $\bigl\|\widehat{\imse}_m(\bft) - \imse(\bft)\bigr\|_{L^q} \lesssim N^{-2\nu/d}$, then
    \begin{align}
        \bigl\|\widehat{\imse}_m(\bft) \bigr\|_{L^q}
        \leq 
        \bigl\| \imse(\bft)\bigr\|_{L^q} 
        + 
        \bigl\|\widehat{\imse}_m(\bft) - \imse(\bft)\bigr\|_{L^q} 
        \lesssim N^{-2\nu/d}. \nonumber
    \end{align}
    In another word, $\widehat{\imse}_m(\bft)$ will achieve the same convergence rate as $\imse(\bft)$.
    Taking $N^{-2\nu/d}$ as a benchmark convergence rate, one can attain this rate by setting
    \begin{align*}
        L &=  \sqrt{2d\nu}\ell \log N, 
        \quad
        m \gtrsim N^{\frac{1}{d} - \frac{1}{2q\nu}} \log N. 
    \end{align*}
    For $q \in [1, \infty]$, the required number of HSGP basis $m$ increases with $q$. 
\end{remark}

\subsection{Error rates for HSGP-based sequential designs with noise/nugget} \label{subsec:nugget}
In practice, the observation $y$ may be noisy, or a nugget term may be added to the kernel $k$ as in \eqref{eq:model_gpr} to ensure numerical stability. In either case, $\eta > 0$ in model \eqref{eq:model_gpr}. The convergence rate of the HSGP approximation for $\imse$ can still be established. In fact, if the nugget parameter $\eta$ decreases sufficiently fast, the same asymptotic rate as in the noiseless case can be recovered. In contrast, if the observation noise variance is bounded away from zero, the asymptotic behavior of IMSE becomes qualitatively different.

If $\eta > 0$, according to Theorem~\ref{thm:imse_acq}, the $\imse$ acquisition takes the form
\begin{align}
    \imse_\eta(\bft)
    &:= \frac{1}{P_{N, \eta}^2(\bft) + \eta}
    \int_{\Omega}\left[k(\bfx,\bft)-\bfk_N^\T(\bfx) (\K_N + \eta \I_N)^{-1} \bfk_N(\bft)\right]^2  \ud \PP(\bfx) \nonumber\\
    &= \frac{1}{P_{N, \eta}^2(\bft) + \eta}
    \int_{\Omega} \left[k(\bfx,\bft) - \Pi_{\X_N, \eta} k(\bfx,\bft)\right]^2 \ud \PP(\bfx), \nonumber
\end{align}
where $[\Pi_{\X_N, \eta} f](\cdot) := \bfk_N^\T(\cdot) (\K_N + \eta \I_N)^{-1} \left[f(\bfx_1), \ldots, f(\bfx_N)\right]^\T$. In this case, the numerator admits the following convergence rate.

\begin{theorem} \label{thm:error_imse_hsgp_numerator_nugget}
    With the same assumptions and notations in Theorem~\ref{thm:error_imse_hsgp_numerator}, for any $\eta > 0$,
    \begin{align}
        &~\quad \left\| \int_\Omega \left\{[\hat{k}_m(\bfx,\cdot) - \Pi_{\X_N, \eta} \hat{k}_m(\bfx,\cdot)]^2 - [k(\bfx,\cdot) - \Pi_{\X_N, \eta} k(\bfx,\cdot)]^2 \right\} \ud \PP(\bfx) \right\|_{L^q} \nonumber\\
        &\lesssim h_N^{d/q} \left(h_N^{\nu} + \eta^{1/2}\right)^2 \left[ \exp \left(-\frac{\sqrt{2\nu}}{2\sqrt{d} \ell} L\right) + (L/m)^{\nu} \right]^2 \nonumber\\
        &\qquad+ h_N^{\frac{d}{2q}} \left(h_N^{\nu} + \eta^{1/2}\right)^2 \min \left\{h_N^{\frac{d}{2q}}, \left(h_N^{\nu} + \eta^{1/2}\right) \right\} \left[\exp \left(-\frac{\sqrt{2\nu}}{2\sqrt{d} \ell} L\right) + (L/m)^{\nu} \right]. \nonumber
    \end{align}
\end{theorem}

This theorem extends Theorem~\ref{thm:error_imse_hsgp_numerator} to the setting $\eta>0$. 
It shows that, after replacing the exact kernel $k$ by its HSGP approximation $\hat{k}_m$, the resulting error in the integrated squared residual term remains controlled by two factors: the interpolation accuracy under the regularized projector $\Pi_{\X_N,\eta}$ and the kernel approximation error.
Compared with the noise-free case, the nugget $\eta$ weakens the interpolation accuracy through the additional $\eta^{1/2}$ term, but the overall structure of the bound remains the same: the error still decays as the design becomes denser ($h_N\to 0$) and as the HSGP approximation improves (by increasing $L$ and $m$).
For the denominator term $P_{N, \eta}^2(\bft) + \eta$, we have the following analysis.
\begin{theorem} \label{thm:lower_bound_power_with_nugget}
    With the same conditions in Theorem~\ref{thm:lower_bound_power_function}, for any $\bft \in \Omega_\gamma$,
    \begin{align}
        P_{N, \eta}^2(\bft) \gtrsim h_N^{2\nu}. \nonumber
    \end{align}
\end{theorem}

By Theorem~\ref{thm:lower_bound_power_with_nugget}, on the restricted set $\Omega_\gamma$, one has $P_{N, \eta}^2(\bft) + \eta \gtrsim \max\{\eta, h_N^{2\nu}\}$ for Mat\'ern kernels with fixed smoothness $\nu$. Then the lower bound of denominator is determined by the order of $h_N^{2\nu}$ and $\eta$. Combining Theorem~\ref{thm:error_imse_hsgp_numerator_nugget} and Theorem~\ref{thm:lower_bound_power_with_nugget}, the upper bound of $\imse$ with regularization follows.
\begin{corollary} [Upper bound for $\widehat{\imse}_m(\bft)$ with noise/nugget case]\label{cor:error_imse_hsgp_numerator_nugget}
    Suppose conditions in Theorems~\ref{thm:error_imse_hsgp_numerator_nugget} and Theorem~\ref{thm:lower_bound_power_with_nugget} hold. With same notations, the followings hold.
    If $\eta \lesssim h_N^{2\nu}$, we obtain results consistent with Corollary~\ref{cor:global_bound_imse_hsgp},
    \begin{align}
        \bigl\|\widehat{\imse}_{\eta, m}(\bft) - \imse_\eta(\bft)\bigr\|_{L^q}
        &\lesssim
        h_N^{\frac{d}{q}} \left[ \exp \left(-\frac{\sqrt{2\nu}}{2\sqrt{d} \ell} L\right) + (L/m)^{\nu} \right]^2  \nonumber\\
        &\quad + h_N^{\nu + \frac{d}{2q}} \left[\exp \left(-\frac{\sqrt{2\nu}}{2\sqrt{d} \ell} L\right) + (L/m)^{\nu} \right]. \nonumber
    \end{align}
    If $\eta \gtrsim h_N^{2\nu}$,
    \begin{align}
        \bigl\|\widehat{\imse}_{\eta, m}(\bft) - \imse_\eta(\bft)\bigr\|_{L^q}
        &\lesssim
        h_N^{\frac{d}{q}} \left[ \exp \left(-\frac{\sqrt{2\nu}}{2\sqrt{d} \ell} L\right) + (L/m)^{\nu} \right]^2  \nonumber\\
        &\quad + h_N^{\frac{d}{2q}} \min \left\{h_N^{\frac{d}{2q}}, \eta^{1/2} \right\} \left[\exp \left(-\frac{\sqrt{2\nu}}{2\sqrt{d} \ell} L\right) + (L/m)^{\nu} \right]. \nonumber
    \end{align}
\end{corollary}

\begin{remark}
    Corollary~\ref{cor:error_imse_hsgp_numerator_nugget} enlightens the selection of nugget parameter $\eta$. Suppose $h_N \asymp N^{-1/d}$ and $\eta \asymp N^{-a}$ for some $a \geq 0$, the optimal convergence rate depends on the value of $a$. 
    If $a \geq 2\nu/d$, then $\eta \lesssim h_N^{2\nu}$, then result in Remark~\ref{rem:choice_L_m} is recovered, and 
    \begin{align*}
        \bigl\|\widehat{\imse}_{\eta, m}(\bft) - \imse_\eta(\bft)\bigr\|_{L^q}
        &\lesssim
        N^{-\frac{1}{q}} \left(\frac{\log m}{m}\right)^{2\nu} 
        +
        N^{-\frac{1}{2q} - \frac{\nu}{d}} \left(\frac{\log m}{m}\right)^\nu.
    \end{align*}
    If $a \in (1/q, 2\nu/d)$, then $\eta \gtrsim h_N^{2\nu}$ and $\eta^{1/2} \lesssim h_N^{\frac{d}{2q}}$, thus
    \begin{align*}
        \bigl\|\widehat{\imse}_{\eta, m}(\bft) - \imse_\eta(\bft)\bigr\|_{L^q}
        &\lesssim
        N^{-\frac{1}{q}} \left(\frac{\log m}{m}\right)^{2\nu}
        +
        N^{-\frac{1}{2q} - \frac{a}{2}} \left(\frac{\log m}{m}\right)^{\nu}.
    \end{align*}
    If $a \in [0, 1/q]$, for example, $\eta \asymp 1$, i.e., $\eta$ is a positive constant, then the aliasing error term is absorbed. Consequently, by letting $L= \sqrt{2d\nu}\ell \log m$, 
    \begin{align*}
        \bigl\|\widehat{\imse}_{\eta, m}(\bft) - \imse_\eta(\bft)\bigr\|_{L^q}
        &\lesssim 
        h_N^{\frac{d}{q}} \left[ \exp \left(-\frac{\sqrt{2\nu}}{2\sqrt{d} \ell} L\right) + (L/m)^{\nu} \right] \lesssim N^{-\frac{1}{q}} \left(\frac{\log m}{m}\right)^\nu.
    \end{align*}
\end{remark}

\vspace{8mm}

\section{Numerical experiments} \label{sec:experiments}
We evaluate the proposed HSGP approximation from two perspectives. First, the approximation accuracy of $\widehat{\imse}_m(\bft)$ is assessed directly by comparing it to a numerically integrated reference $\imse(\bft)$ over candidate grids. Secondly, the practical impact on sequential design is examined by running Algorithm~\ref{alg:sequential_with_HSGP-IMSE} and checking whether the HSGP surrogate preserves the acquisition landscape that drives point selection.

\subsection{HSGP approximation of IMSE} \label{subsec:experiment.imse}
The first experiment studies the pointwise agreement between the exact acquisition $\imse(\bft)$ and its HSGP approximation $\widehat{\imse}_m(\bft)$. The reference $\imse(\bft)$ is computed via numerical integration on $\Omega$, whereas $\widehat{\imse}_m(\bft)$ is computed using the closed-form expression in Theorem~\ref{thm:HSGP-IMSE_acq}. The comparison is repeated for several commonly used covariance kernels, covering Gaussian processes with both infinite smoothness (Gaussian) and finite-order smoothness (Mat\'ern and generalized Wendland).

\begin{figure}[t]
    \centering
    \begin{subfigure}[t]{0.32\linewidth}
        \centering
        \includegraphics[width=\linewidth]{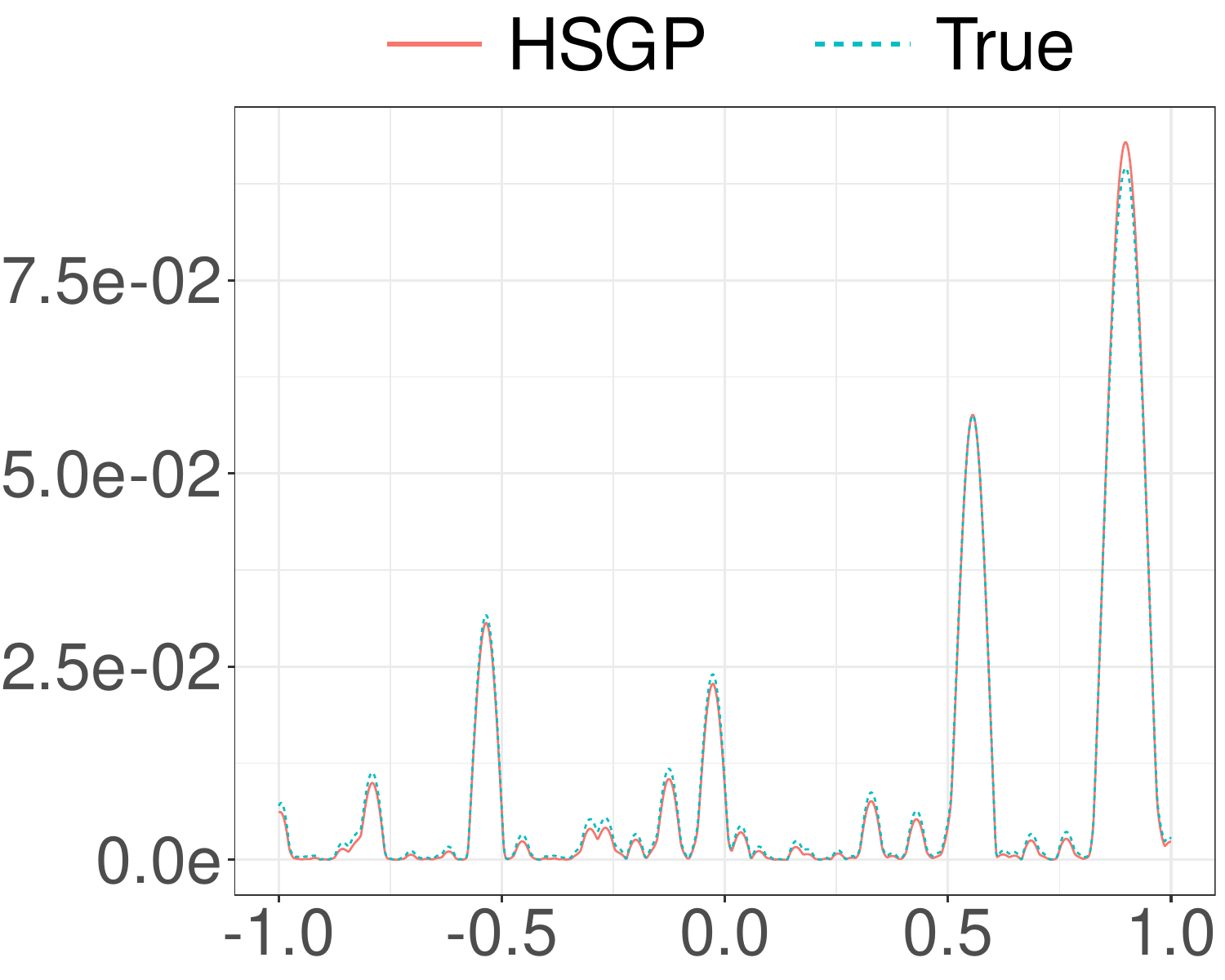}
        \caption{Generalized Wendland (1D).}
        \label{fig:imse_1d_gw}
    \end{subfigure}\hfill
    \begin{subfigure}[t]{0.32\linewidth}
        \centering
        \includegraphics[width=\linewidth]{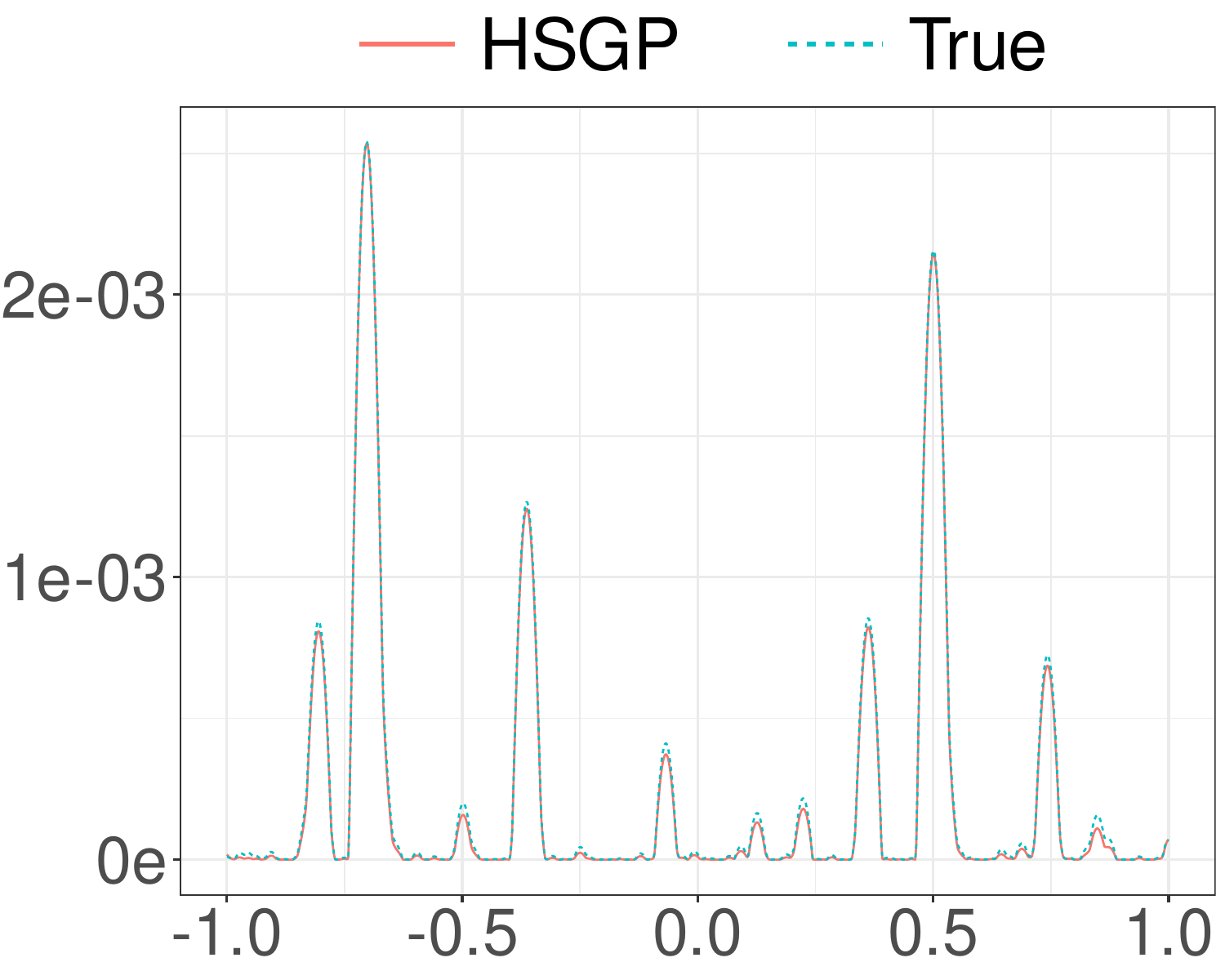}
        \caption{Mat\'ern-$3/2$ (1D).}
        \label{fig:imse_1d_matern}
    \end{subfigure}\hfill
    \begin{subfigure}[t]{0.32\linewidth}
        \centering
        \includegraphics[width=\linewidth]{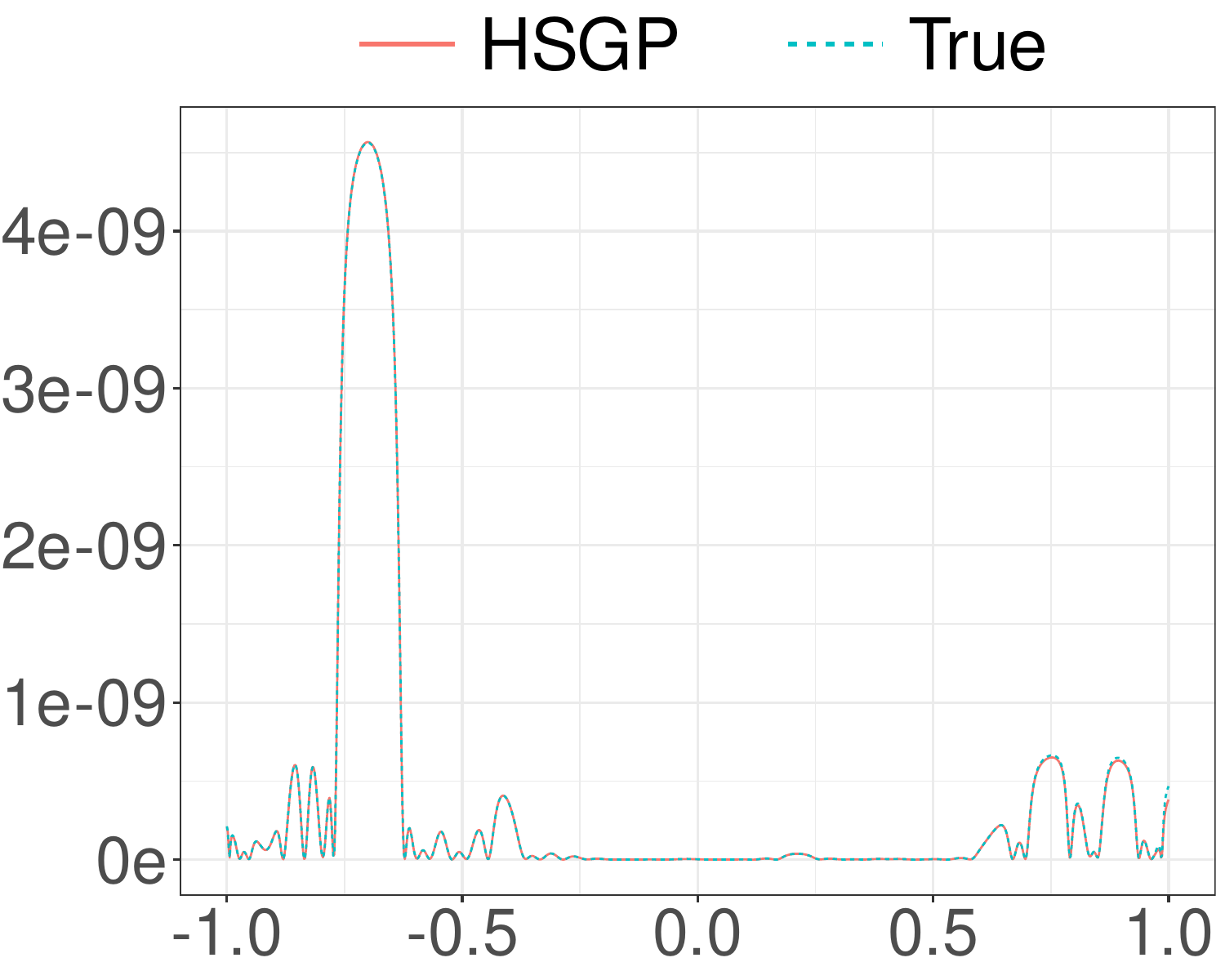}
        \caption{Gaussian (1D).}
        \label{fig:imse_1d_gaussian}
    \end{subfigure}
    \caption{Comparison of the numerically integrated reference acquisition $\imse(\bft)$ and the HSGP approximation $\widehat{\imse}_m(\bft)$ over a one-dimensional candidate grid. Agreement of the two curves indicates that the HSGP surrogate accurately reproduces the acquisition profile used for selecting new design points.}
    \label{fig:imse_1d_three_kernels}
\end{figure}

\paragraph{One-dimensional domain.} Figure~\ref{fig:imse_1d_three_kernels} plots $\imse(\bft)$ and $\widehat{\imse}_m(\bft)$ over a one-dimensional candidate grid for three kernels. The fixed domain $\Omega = (-1, 1)$, i.e., $d=B=1$. The samples used to draw plots are from the Latin Hypercube design, a space filling design. The scaling parameter $\sigma^2$ is set as a constant 2. For numerical stability, the nugget parameter $\eta$ is set as fixed constant $\sigma^2 \cdot 10^{-10}$ to guarantee the positiveness of covariance matrix. 
The Mat\'ern-$3/2$ kernel is assigned length scale $\ell = 0.1$, $N=200$ samples and HSGP parameters $m = 120, L = 1.5$. The generalized Wendland kernel is tested with $\mu = 7, \kappa = 1, \beta = 0.2$, $N=100$ and HSGP parameters $m = 100, L = 1.2$. The parameterization of generalized Wendland kernel is based on Equation (3) in \citet{bevilacqua:2019}. The implementation is through the R package \texttt{GeneralizedWendland} (\citealt{fischer:2022}). The Gaussian kernel is tested with length scale $\ell = 0.1$, $N=100$ and HSGP parameters $m = 100, L = 2$. In this experiment, the two curves of $\imse(\bft)$ and $\widehat{\imse}_m(\bft)$ almost completely overlap across all the three tested cases, showing that the HSGP approximation is highly accurate. 

\begin{figure}[!t]
    \centering

    \begin{subfigure}[t]{0.32\linewidth}
        \centering
        \includegraphics[width=\linewidth]{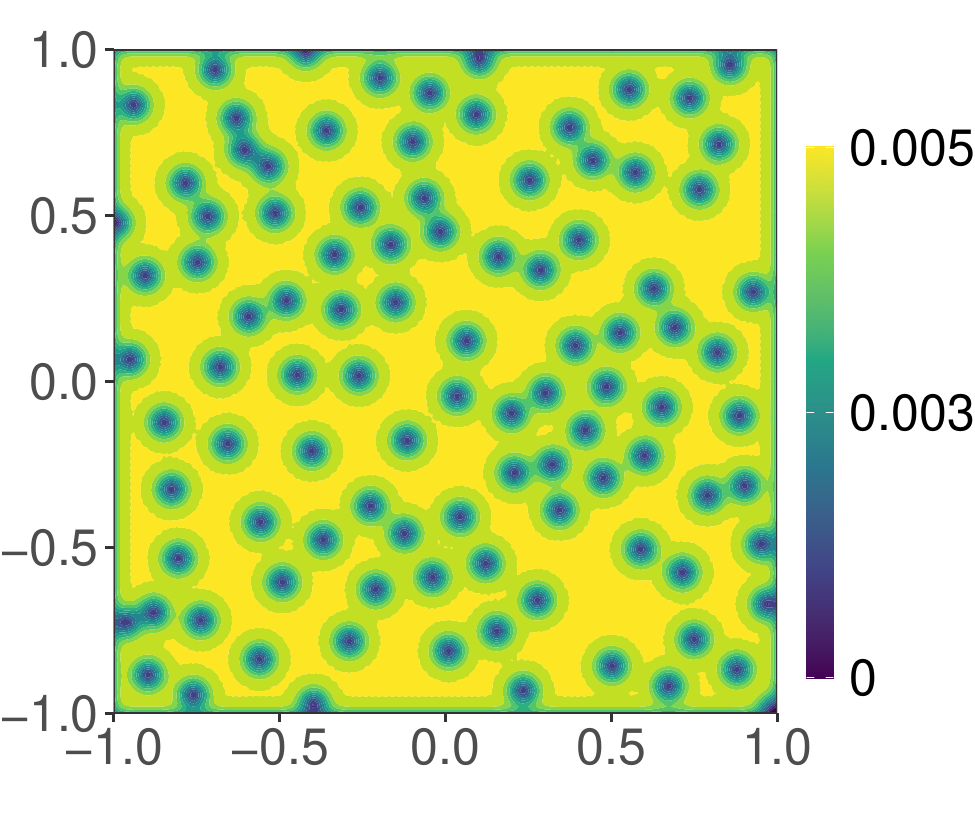}
        \caption{$\imse$ (GW).}
        \label{fig:imse2d_true_gw}
    \end{subfigure}\hfill
    \begin{subfigure}[t]{0.32\linewidth}
        \centering
        \includegraphics[width=\linewidth]{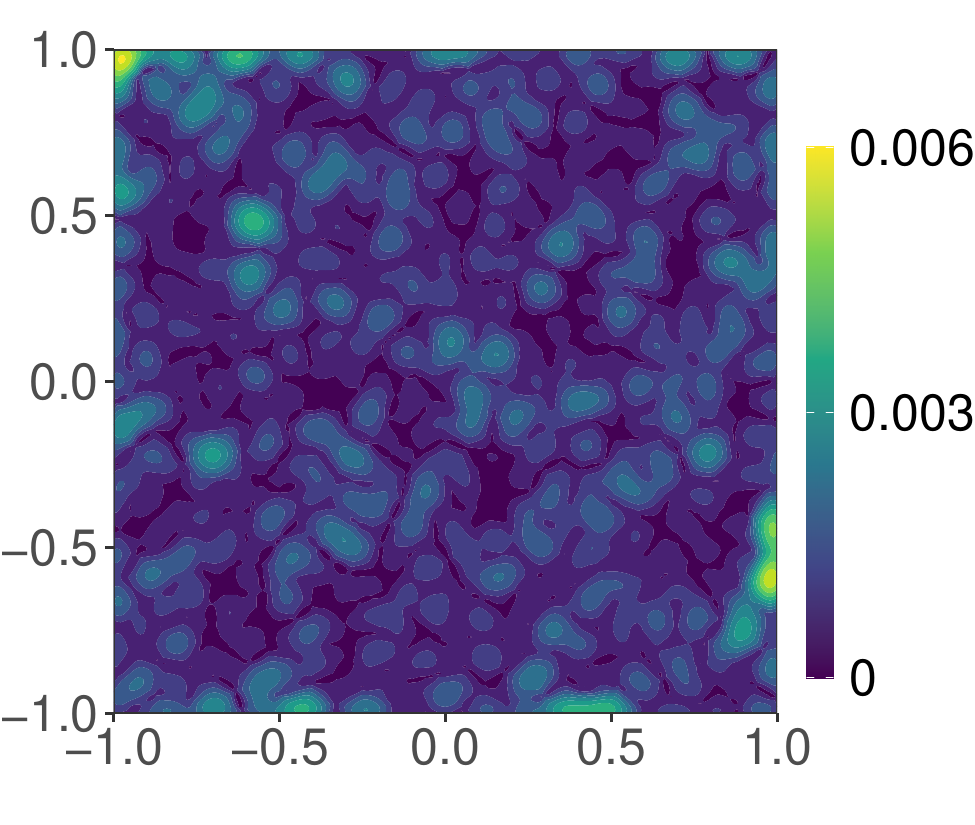}
        \caption{$\imse$ (Mat\'ern-$3/2$).}
        \label{fig:imse2d_true_matern}
    \end{subfigure}\hfill
    \begin{subfigure}[t]{0.32\linewidth}
        \centering
        \includegraphics[width=\linewidth]{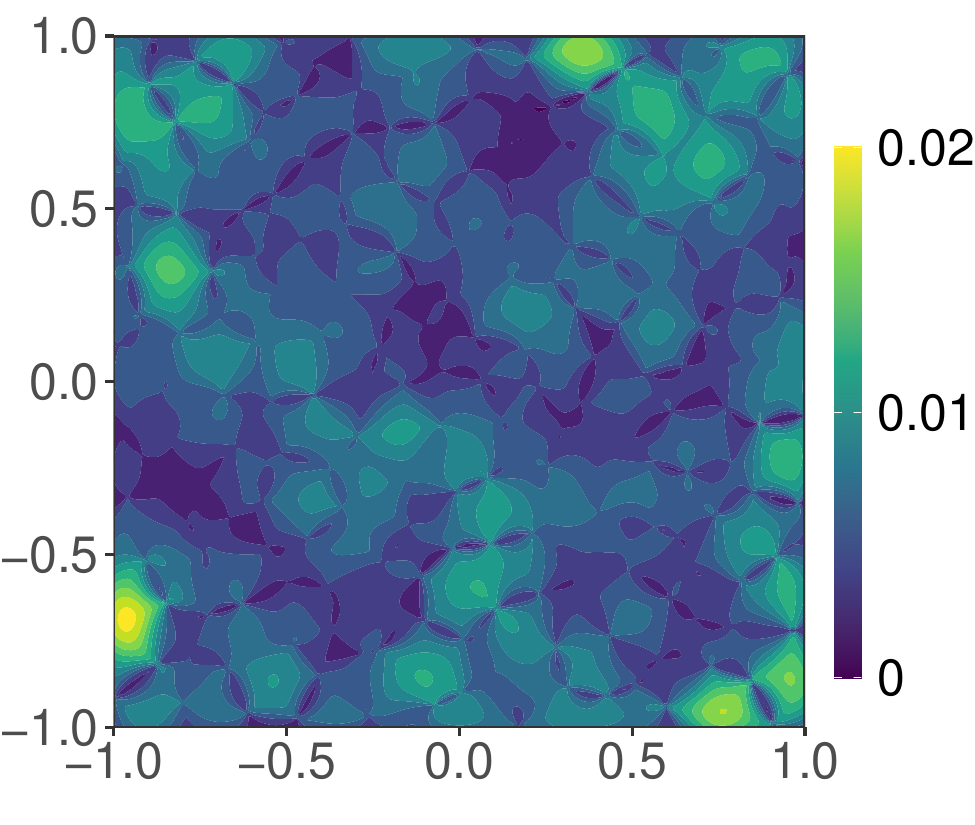}
        \caption{$\imse$ (Gaussian).}
        \label{fig:imse2d_true_gauss}
    \end{subfigure}

    \vspace{0.6em}

    \begin{subfigure}[t]{0.32\linewidth}
        \centering
        \includegraphics[width=\linewidth]{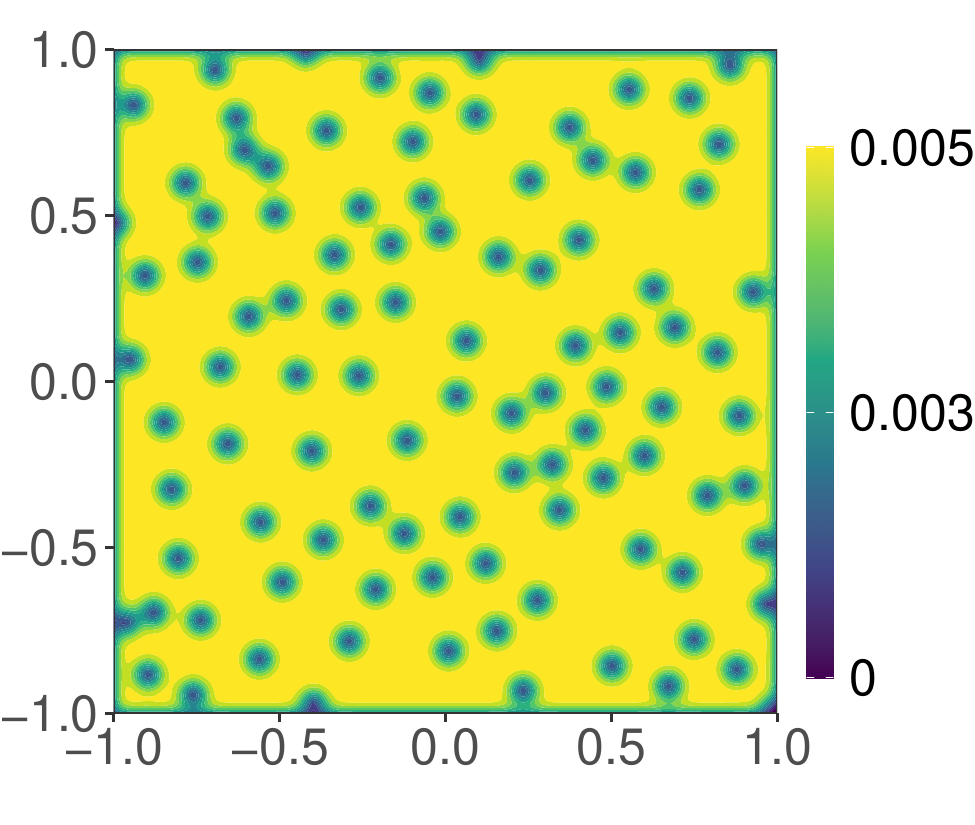}
        \caption{$\widehat{\imse}_m$ (GW).}
        \label{fig:imse2d_hsgp_gw}
    \end{subfigure}\hfill
    \begin{subfigure}[t]{0.32\linewidth}
        \centering
        \includegraphics[width=\linewidth]{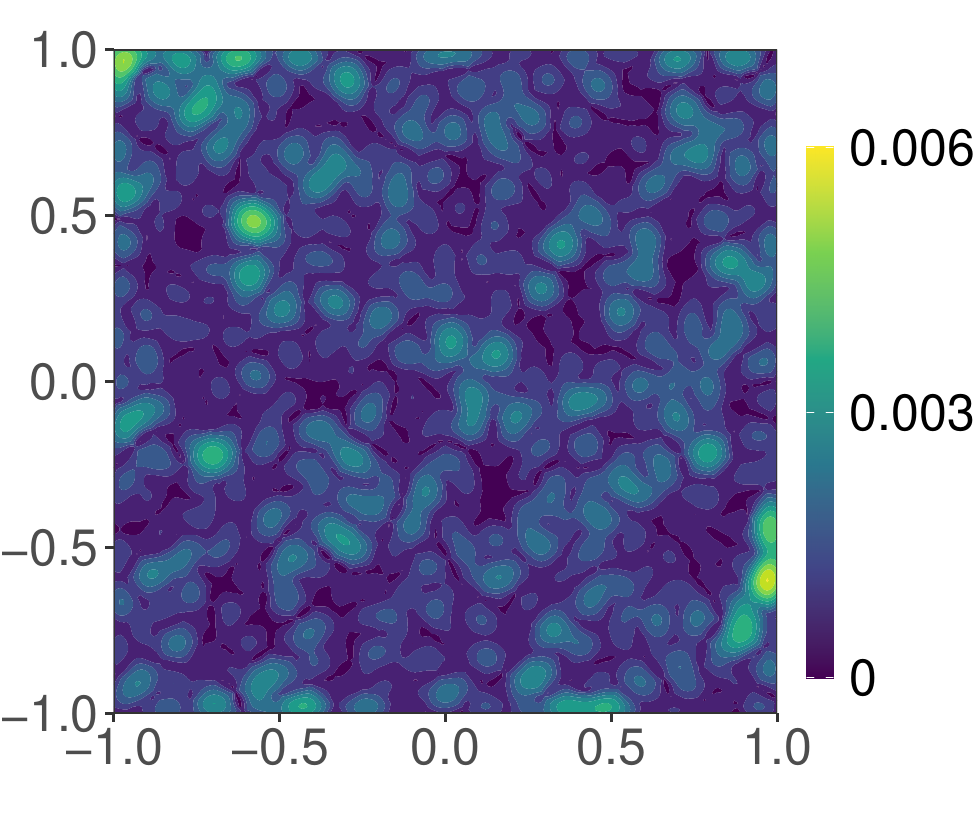}
        \caption{$\widehat{\imse}_m$ (Mat\'ern-$3/2$).}
        \label{fig:imse2d_hsgp_matern}
    \end{subfigure}\hfill
    \begin{subfigure}[t]{0.32\linewidth}
        \centering
        \includegraphics[width=\linewidth]{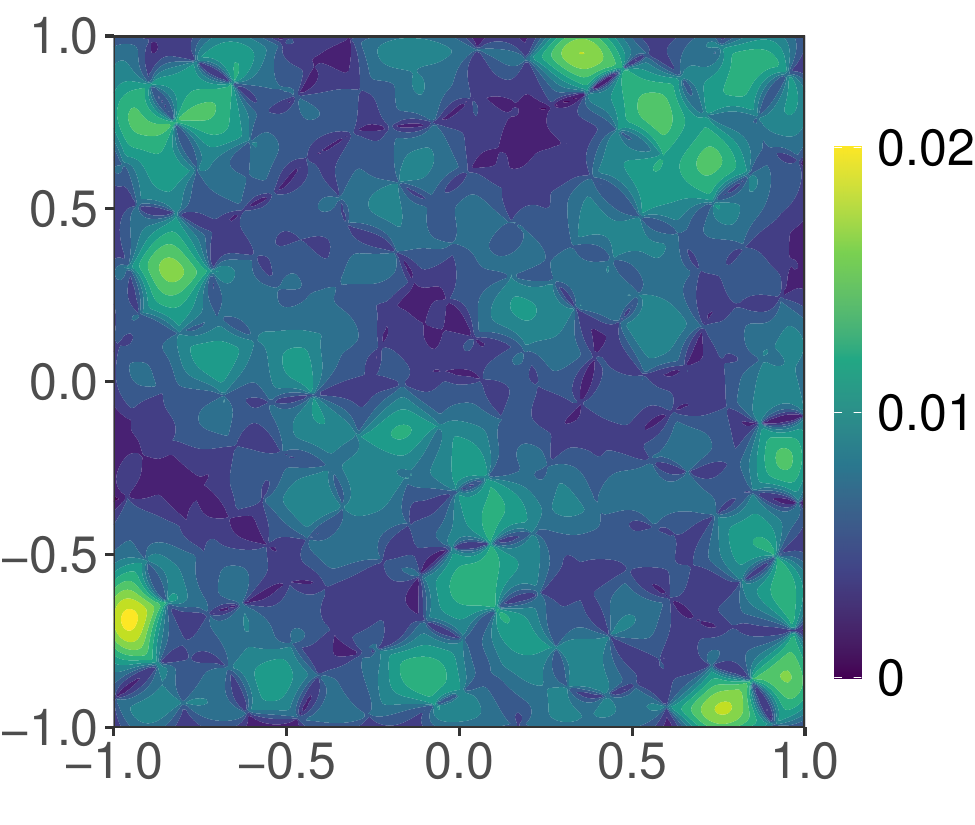}
        \caption{$\widehat{\imse}_m$ (Gaussian).}
        \label{fig:imse2d_hsgp_gauss}
    \end{subfigure}

    \caption{Two-dimensional acquisition contours on $\Omega=(-1,1)^2$ for three covariance kernels. The top row shows the numerically integrated reference acquisition $\mathrm{IMSE}(\bft)$, and the bottom row shows the HSGP approximation $\widehat{\imse}_m(\bft)$ computed from Theorem~\ref{thm:HSGP-IMSE_acq}. Consistency between corresponding panels indicates that the HSGP approximation preserves the acquisition landscape that guides sequential sampling, including the location and geometry of high-acquisition regions.}
    \label{fig:imse2d_true_vs_hsgp_2x3}
\end{figure}

\paragraph{Two-dimensional domain.} For $d=2$, the acquisition function becomes a surface over $\Omega=(-1,1)^2$, and the geometry of level sets is critical: sequential design depends on whether maxima, ridges, and basins are preserved. Figure~\ref{fig:imse2d_true_vs_hsgp_2x3} compares contour plots of $\imse(\bft)$ (top row) against $\widehat{\imse}_m(\bft)$ (bottom row) for the three kernels. The other parameter settings remain the same as in the one-dimensional case. The HSGP approximation is considered effective when the bottom row reproduces the dominant geometric features of the top row, including the location of prominent modes and the overall spread of high-acquisition regions.

In summary, Figures~\ref{fig:imse_1d_three_kernels} and \ref{fig:imse2d_true_vs_hsgp_2x3} are designed to validate two practical requirements for HSGP-based sequential design:
(i) \emph{pointwise accuracy} of $\widehat{\imse}_m(\bft)$ as a surrogate for $\imse(\bft)$, and
(ii) \emph{geometric fidelity} of the acquisition surface, so that optimizing the surrogate yields similar maximizers to optimizing the reference acquisition. These requirements motivate the subsequent experiments on end-to-end sequential simulation performance.

\subsection{Sequential simulations} \label{subsec:experiment.simulation}
This subsection evaluates end-to-end sequential design performance in simulation examples with one-dimensional and two-dimensional domains. Performance is summarized over sequential iterations using three criteria: (i) predictive accuracy measured by the test-set \emph{root mean square prediction error} (RMSE), (ii) uncertainty quantification summarized by the mean posterior variance over test inputs, and (iii) cumulative time cost. The central objective is to assess whether replacing the exact acquisition $\imse(\bft)$ with its HSGP surrogate $\widehat{\imse}_m(\bft)$ preserves the design quality while improving computational efficiency.

In the sequential experiment, the kernel parameters $\sigma^2$ and $\ell$ are fitted by MLE estimators \eqref{eq:log-likelihood}--\eqref{eq:mle_numerical}. In the one dimensional experiment with noisy observations, the noise is assumed Gaussian with zero mean but the variance is unknown. In this case, the noise parameter $\eta$ will also be reparameterized and fitted by MLE estimator \eqref{eq:log-likelihood}--\eqref{eq:mle_numerical}. While for noiseless cases, we set $\eta = \sigma^2 \cdot 10^{-10}$.

For the Gaussian, Mat\'ern-$3/2$, and Mat\'ern-$5/2$ kernels, $\hsgp$-$\imse$ is compared against the IMSPE method in \citet{binois:2019}, using the implementation in the R package \texttt{hetGP}. For kernel families not supported by the \texttt{hetGP} package, $\hsgp$-$\imse$ is compared against a space-filling baseline, where new points are generated by Latin hypercube sampling (LHS). The HSGP parameters $m = 20 \cdot d + 0.1 \cdot B / \ell \cdot \log(N)$ and $L = B + 0.5 \cdot \ell / B \cdot \log(N)$ are updated within each iteration. Both the parameters $m$ and $L$ increase in the order of $\log(N)$ and are also dependent on the dimension $d$, length scale $\ell$, and boundary $B$. Each configuration is repeated over $10$ independent replicates. For each metric, the envelope band displayed in the plots is the pointwise minimum--maximum range across replicates. 

\paragraph{Benchmark functions.} All experiments are conducted on $\Omega=(-1,1)^d$ with the following test functions:
\begin{itemize}
    \item $f_1(\bfx) = \bfk_N^\T(\bfx)\K_N^{-1}\Y_N$, where $\Y_N=[y_1,\ldots,y_N]^\T\sim\mathcal{GP}(0,k)$ is first generated from a Mat\'ern-$5/2$ Gaussian random field with $\sigma^2=1$ and $\ell=0.1$ using $N=200d$ initial locations. The resulting posterior mean function $f_1$ is then treated as a deterministic benchmark.
    \item $f_2(\bfx)=p(\bfx\mid\bm\mu_1,\varrho^2)+p(\bfx\mid\bm\mu_2,\varrho^2)-p(\bfx\mid\bm\mu_3,\varrho^2)-p(\bfx\mid\bm\mu_4,\varrho^2)$,
    where $p(\cdot\mid\bm\mu,\varrho^2)$ is the density of $\mathcal{N}(\bm\mu,\varrho^2\I_2) $,
    $\bm\mu_1=(0.5,0.5)^\T$, $\bm\mu_2=(-0.5,-0.5)^\T$, $\bm\mu_3=(0.5,-0.5)^\T$, $\bm\mu_4=(-0.5,0.5)^\T$, and $\varrho^2=0.01$.
    \item $f_3(\bfx) = \sum_{k=1}^{100} A_k \exp\left(-\|\bfx-\bfc_k\|_2/s_k\right)$, where $A_k\sim\mathrm{Unif}(-1,1)$, $\bfc_k\sim\mathrm{Unif}((-1, 1)^d)$ and $s_k\sim\mathrm{Unif}(0.01,1)$.
    \item $f_4(\bfx) = \sum_{k=1}^d \cos \left(\frac{10\pi x_k}{1+x_k+5x_k^2}\right)$.
\end{itemize}

\begin{figure}[!t]
  \centering
    \includegraphics[width=.5\textwidth,keepaspectratio]{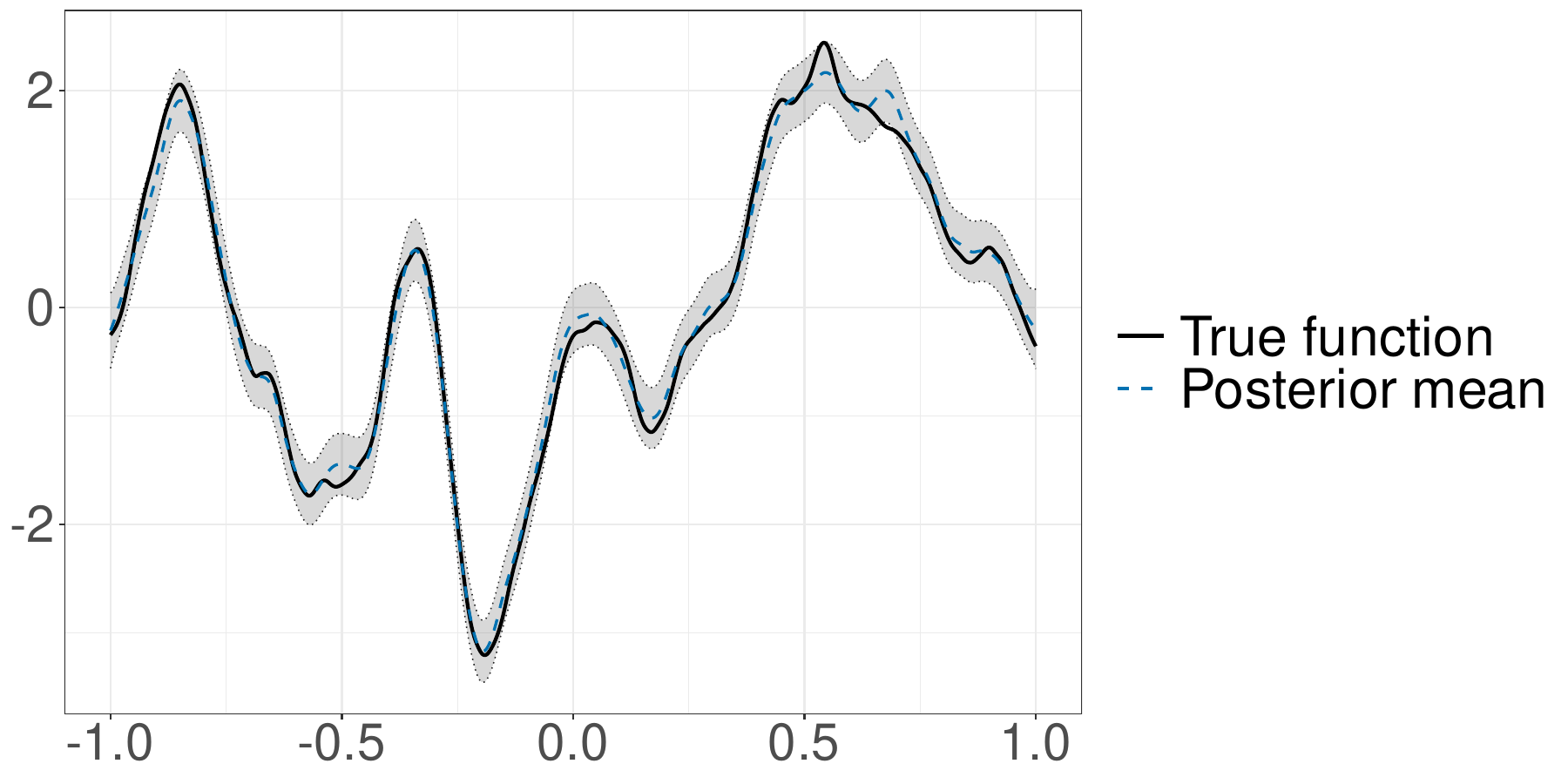}
    \caption{HSGP-IMSE posterior mean vs.\ true function with 95\% confidence interval. This diagnostic plot overlays the fitted GP posterior mean with the ground-truth response on a dense grid.}
    \label{fig:prediction_1d_matern3_2}
\end{figure}

\paragraph{One-dimensional Mat\'ern-$3/2$ kernel.} In the one-dimensional simulation, 100 initial samples are generated from $y = f_1(\bfx) + \epsilon, \epsilon \sim \mathcal{N}(0, 0.025)$. The locations $\bfx$ are sampled from $\Omega=(-1,1)$ using LHS. Then Algorithm~\ref{alg:sequential_with_HSGP-IMSE} is applied to sample 400 more data points using a Mat\'ern-$3/2$ GP surrogate. 

\begin{figure}[!t]
  \centering
  \begin{subfigure}[t]{0.32\textwidth}
    \centering
    \includegraphics[width=\textwidth]{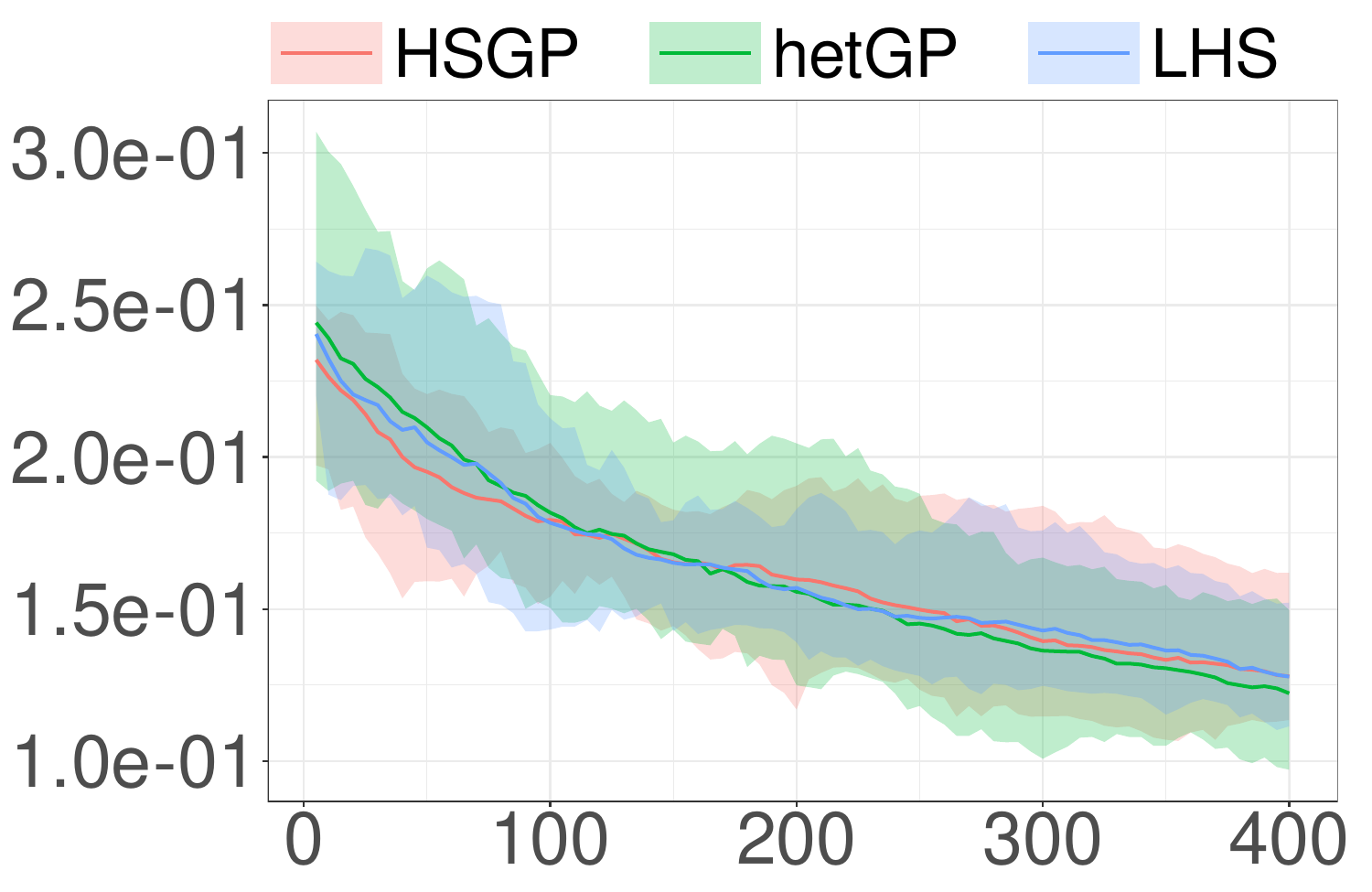}
    \caption{RMSE (Mat\'ern-$3/2$).}
    \label{fig:matern3_2-1d-rmse}
  \end{subfigure}\hfill
  \begin{subfigure}[t]{0.32\textwidth}
    \centering
    \includegraphics[width=\textwidth]{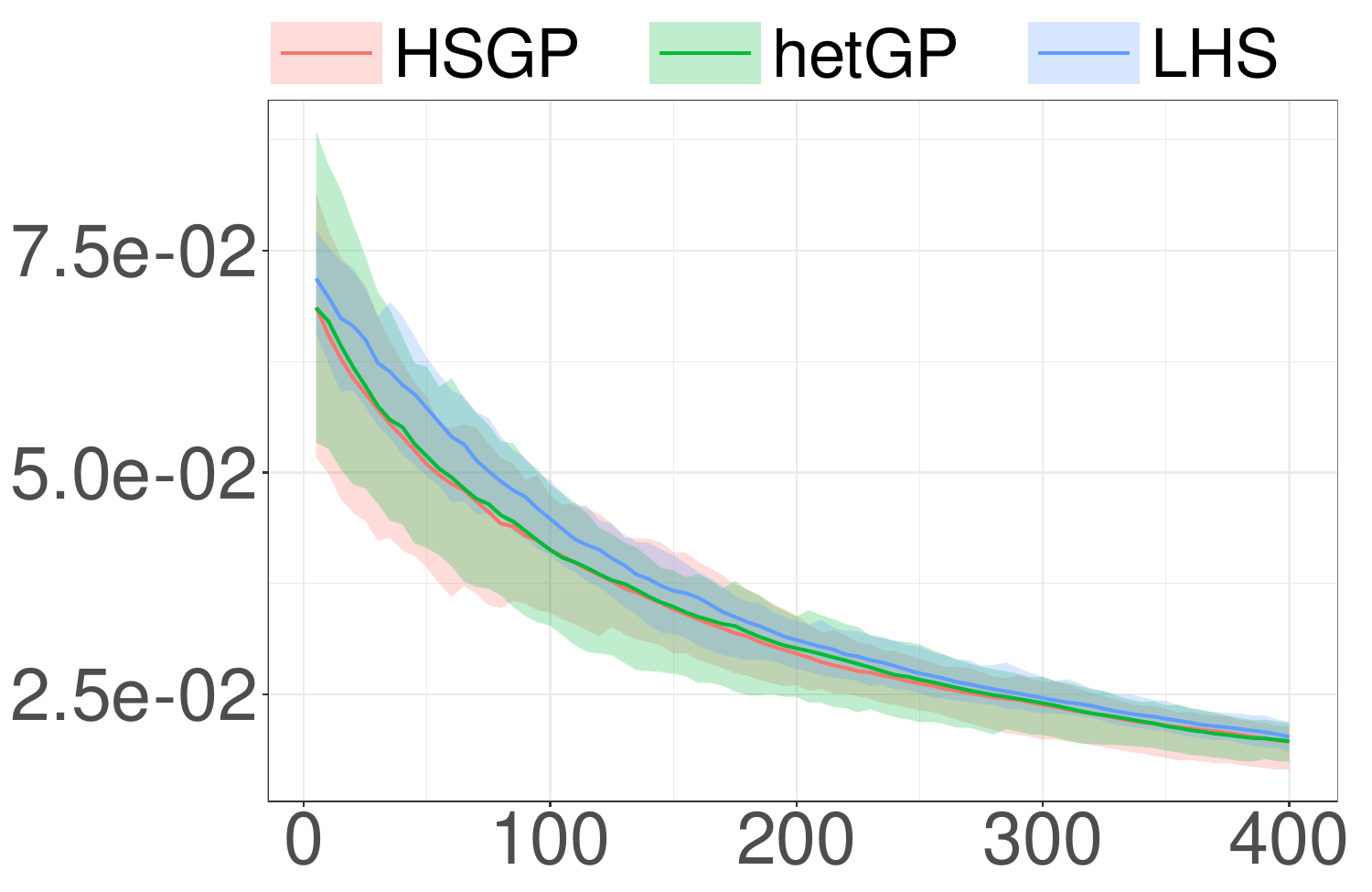}
    \caption{Variance (Mat\'ern-$3/2$).}
    \label{fig:matern3_2-1d-sd2}
  \end{subfigure}\hfill
  \begin{subfigure}[t]{0.32\textwidth}
    \centering
    \includegraphics[width=\textwidth]{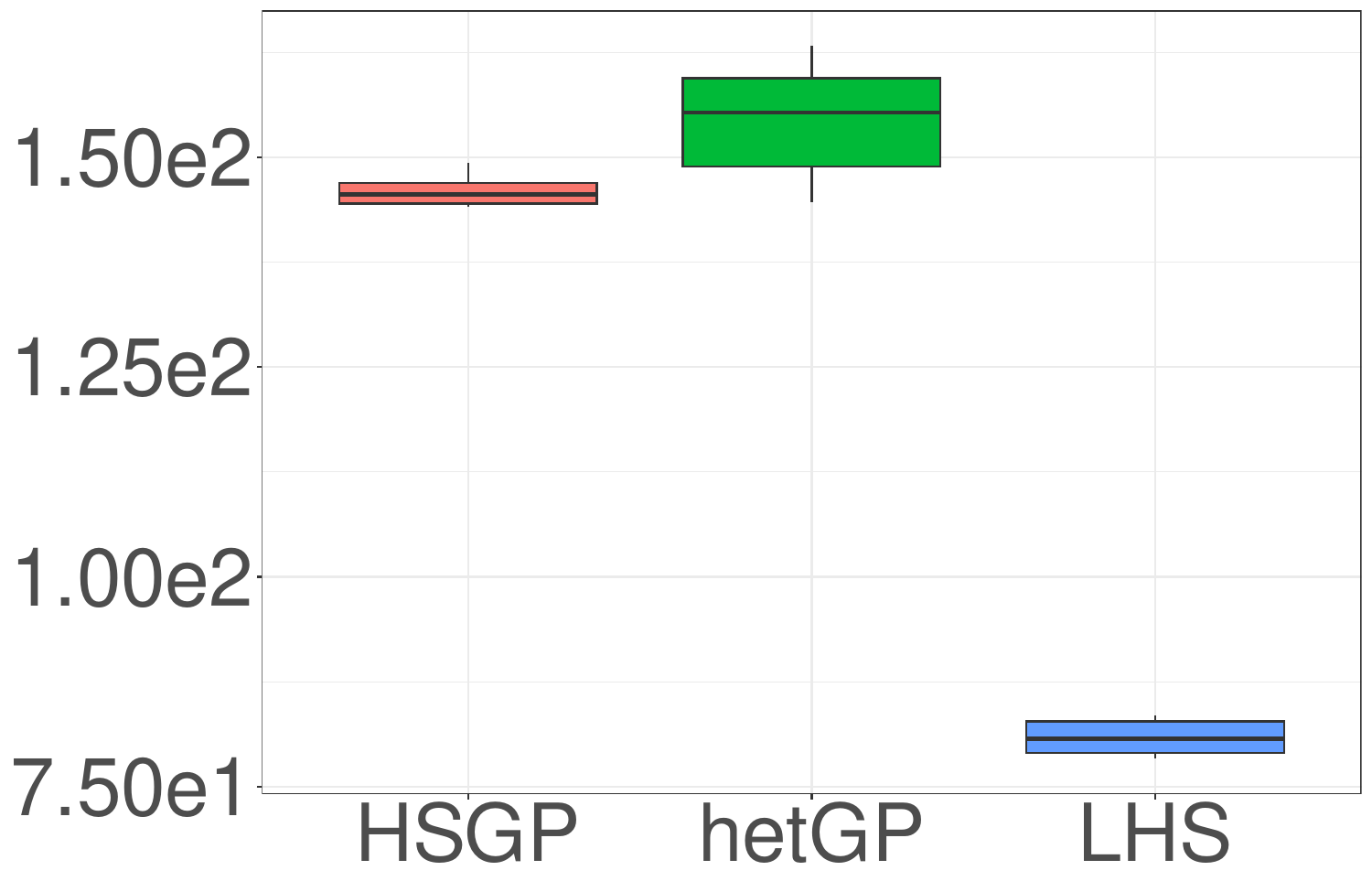}
    \caption{Time cost (Mat\'ern-$3/2$).}
    \label{fig:matern3_2-1d-time}
  \end{subfigure}
  \caption{One-dimensional simulation with a Mat\'ern-$3/2$ kernel. Curves show averaged RMSE (left), mean posterior variance (Variance) over test points (middle), and cumulative time cost (right) across sequential iterations. The envelope band denotes the pointwise minimum--maximum across replicates, capturing run-to-run variability induced by random initialization and stochastic optimization.}
  \label{fig:sim_matern3_2_1d}
\end{figure}

The posterior mean and 95\% confidence interval are plotted in Figure~\ref{fig:prediction_1d_matern3_2}. This experiment serves as a controlled setting to verify that $\widehat{\imse}_m$ yields a valid sequential design trajectory. The same data points are also fitted by the IMSPE method in \texttt{hetGP} package and by the LHS design. Figure~\ref{fig:sim_matern3_2_1d} reports the RMSE, mean posterior variance, and time cost when fitting $f_1$ using the three different methods. RMSE directly measures the primary objective of improved predictive accuracy on held-out test points; the mean posterior variance summarizes uncertainty quantification and reflects how quickly the sequential design reduces global uncertainty over $\Omega$; and cumulative time cost isolates the computational trade-off induced by acquisition optimization. Together, these panels align with the goal of assessing whether the surrogate acquisition $\widehat{\imse}_m$ preserves design quality while offering practical speed advantages. In this test, all the three methods achieve similar performance in RMSE and prediction uncertainty reduction from Figure~\ref{fig:matern3_2-1d-rmse} and \ref{fig:matern3_2-1d-sd2}. On the other hand, our HSGP-IMSE method achieves similar time cost compared with IMSPE method from Figure~\ref{fig:matern3_2-1d-time}.

\paragraph{Two-dimensional Gaussian kernel and Mat\'ern-$5/2$ simulations.} The next set of experiments considers the two-dimensional case ($d=2$) and examines two additional kernels supported by the \texttt{hetGP} package. In this setting, the acquisition landscape is more challenging, and optimization of the acquisition function becomes more computationally demanding. The experiments are initialized with 500 LHS samples and proceed sequentially until a total of 2000 samples is reached.

\begin{figure}[!t]
  \centering
  \begin{subfigure}[t]{0.49\textwidth}
    \centering
    \includegraphics[width=\textwidth]{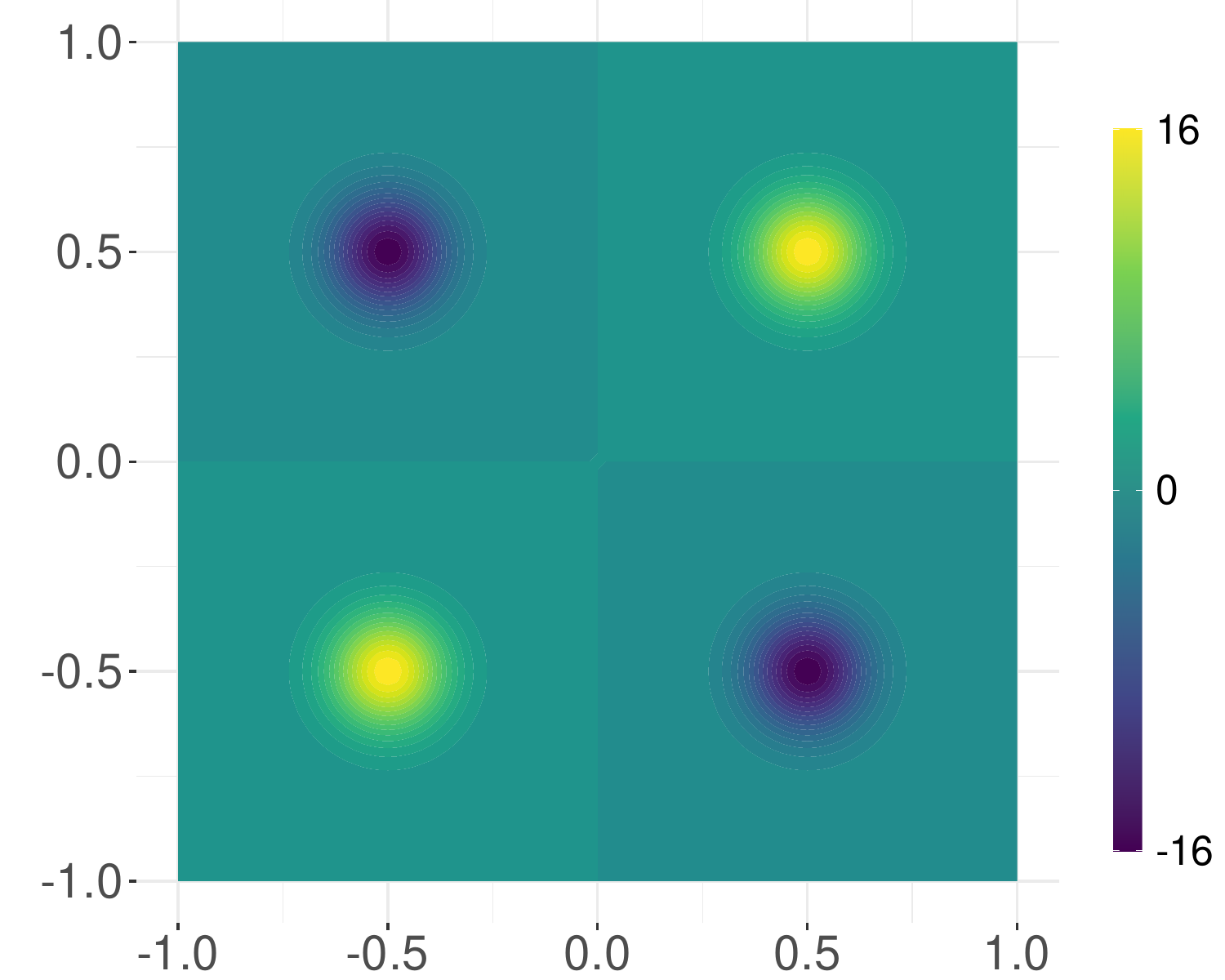}
    \caption{True function $f(\bfx)$ (contour).}
    \label{fig:gaussian-true-contour}
  \end{subfigure}\hfill
  \begin{subfigure}[t]{0.49\textwidth}
    \centering
    \includegraphics[width=\textwidth]{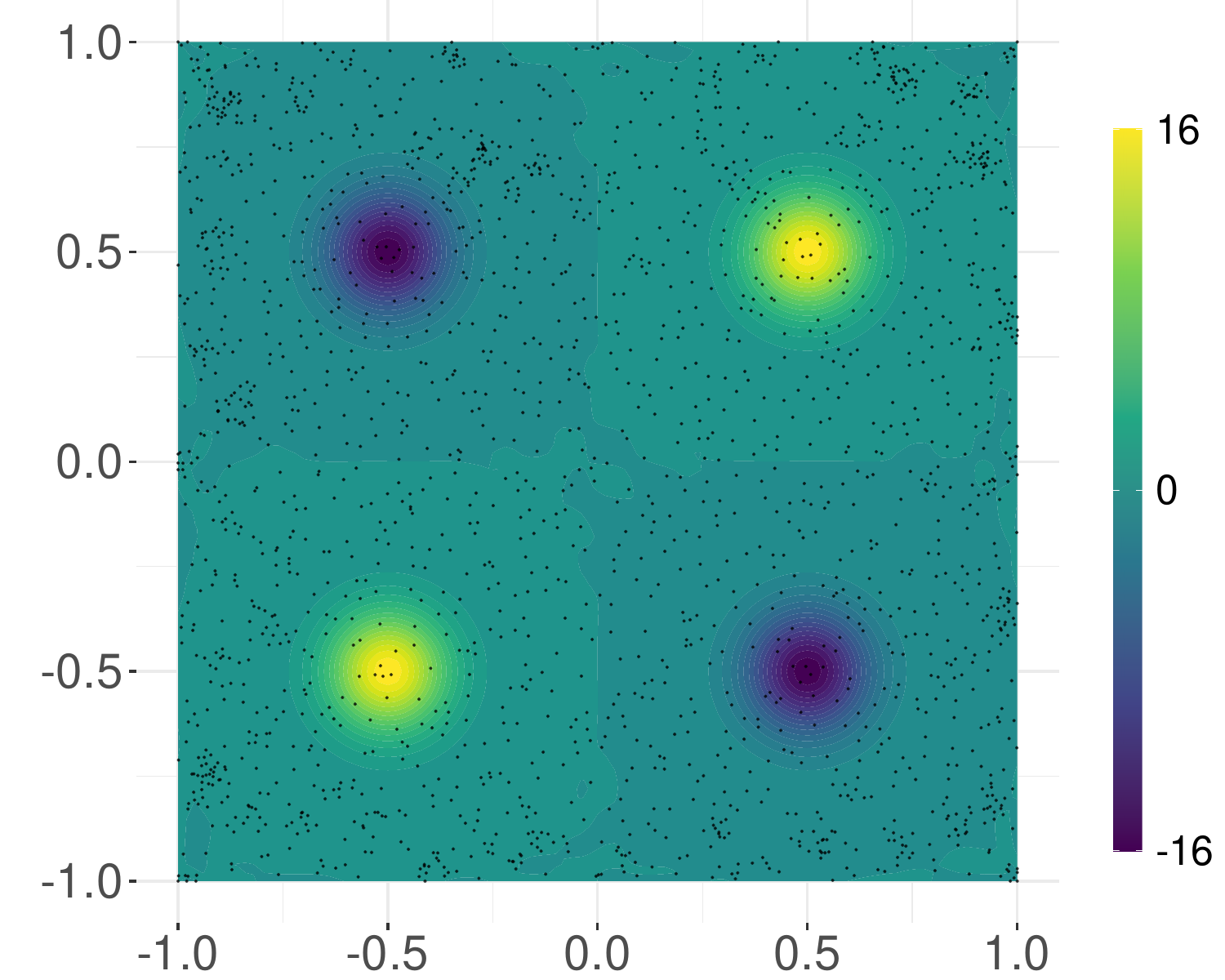}
    \caption{Posterior mean (contour).}
    \label{fig:gaussian-mean-contour}
  \end{subfigure}
  \caption{Two-dimensional Gaussian kernel experiment.}
  \label{fig:prediction3d_gaussian}
\end{figure}

\begin{figure}[!t]
  \centering
  \begin{subfigure}[t]{0.32\textwidth}
    \centering
    \includegraphics[width=\textwidth,]{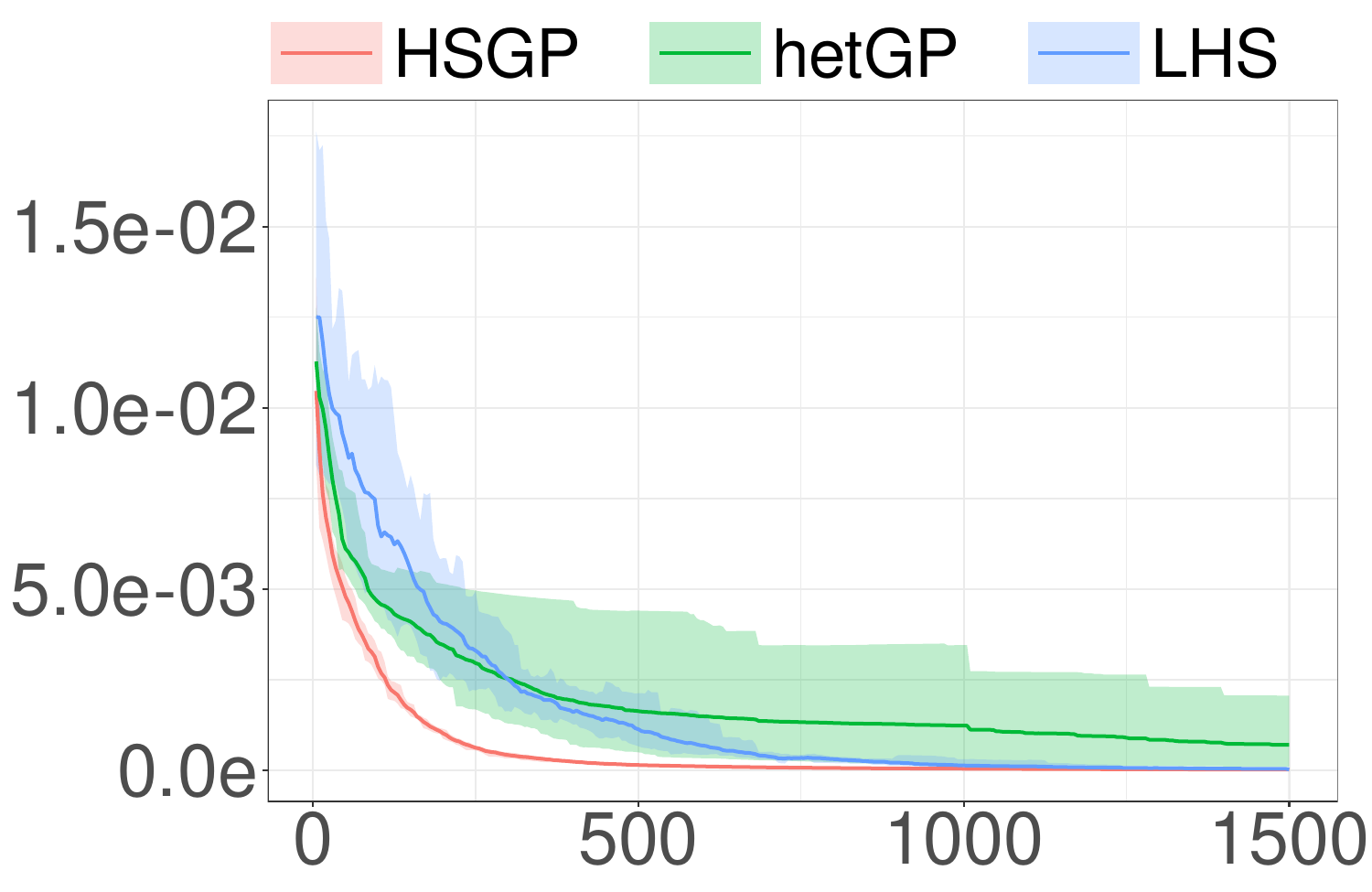}
    \caption{RMSE (Gaussian).}
    \label{fig:gaussian-2d-rmse}
  \end{subfigure}\hfill
  \begin{subfigure}[t]{0.32\textwidth}
    \centering
    \includegraphics[width=\textwidth,]{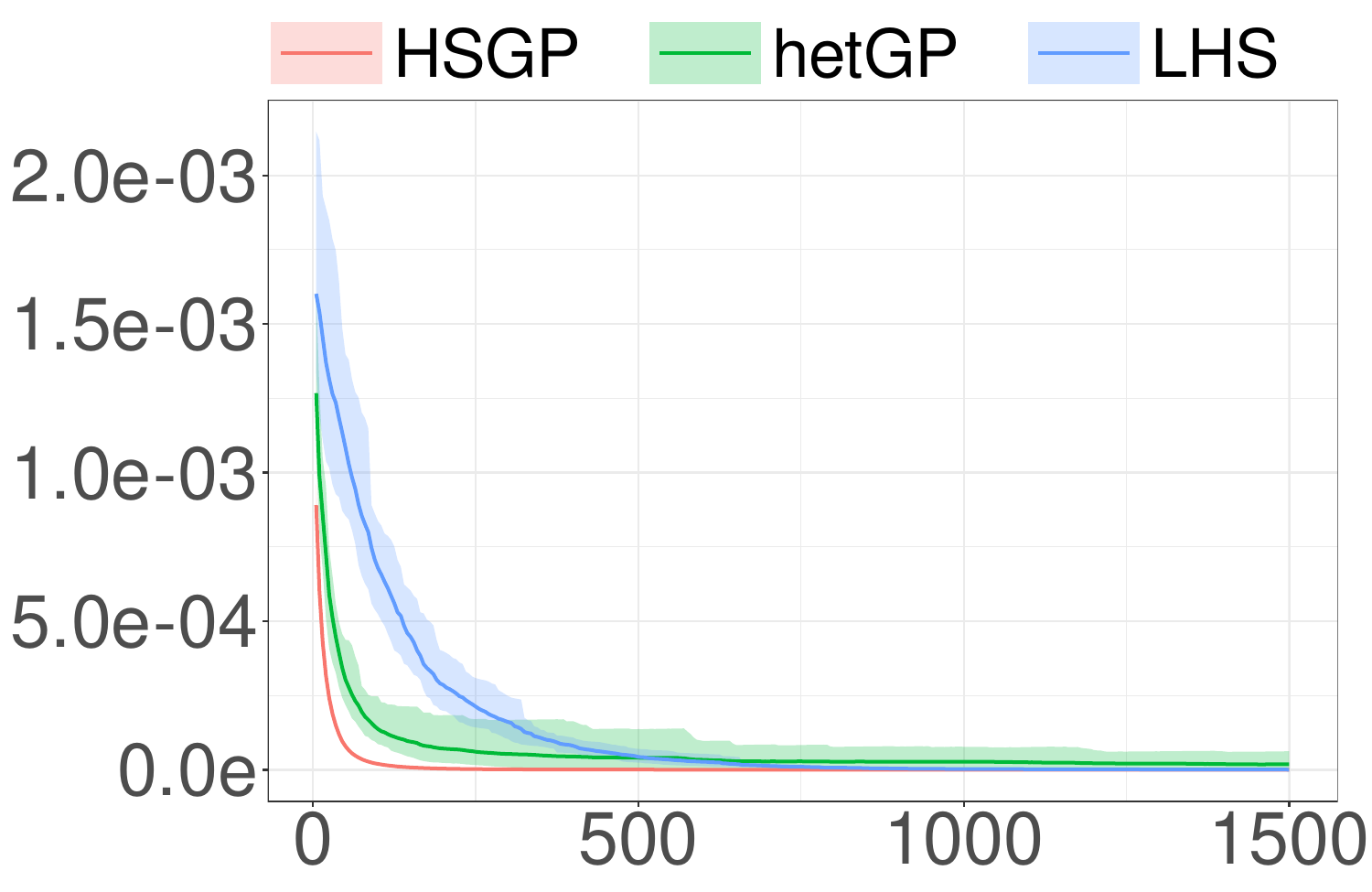}
    \caption{Variance (Gaussian).}
    \label{fig:gaussian-2d-sd2}
  \end{subfigure}\hfill
  \begin{subfigure}[t]{0.32\textwidth}
    \centering
    \includegraphics[width=\textwidth,]{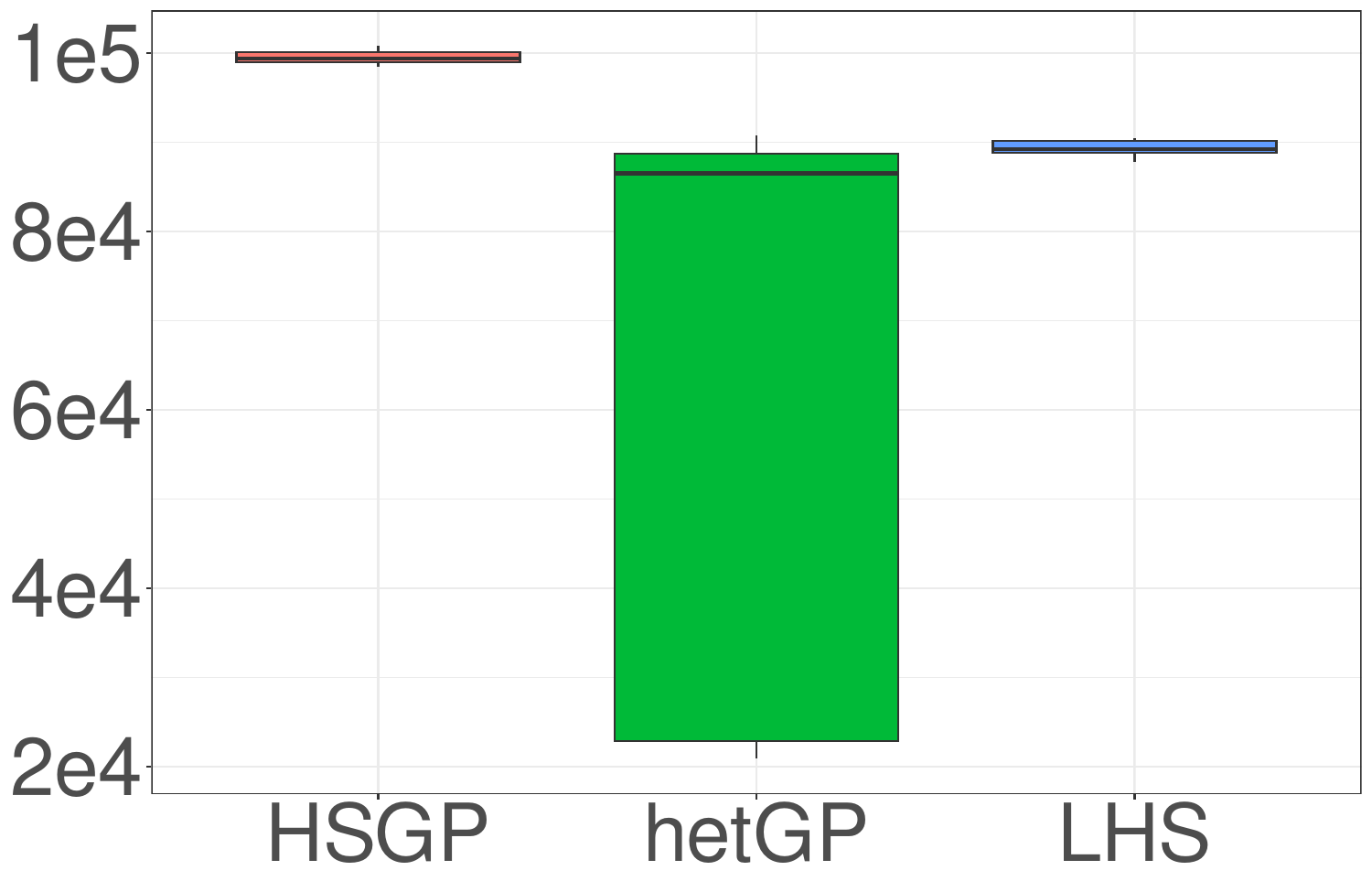}
    \caption{Time cost (Gaussian).}
    \label{fig:gaussian-2d-time}
  \end{subfigure}
  \caption{Two-dimensional simulation using 3 methods with a Gaussian kernel GP surrogate.}
  \label{fig:sim_gaussian_2d}
\end{figure}

First, 500 initial samples are generated from $y = f_2(\bfx)$, and Algorithm~\ref{alg:sequential_with_HSGP-IMSE} is applied using the Gaussian kernel. Figure~\ref{fig:gaussian-true-contour} is the contour plot of the benchmark function $f_2(\bfx)$ on a dense grid. Figure~\ref{fig:gaussian-mean-contour} is the plot of posterior mean after adding 1500 sequentially sampled points using the HSGP-IMSE method. Comparing Figure \ref{fig:gaussian-true-contour} with \ref{fig:gaussian-mean-contour}, it can be concluded that HSGP-IMSE method successfully recovers the true benchmark function. Next, Figure~\ref{fig:sim_gaussian_2d} demonstrates that HSGP-IMSE outperforms IMSPE and LHS in reducing RMSE and prediction uncertainty, while achieving similar time cost.

\begin{figure}[!t]
  \centering
  \begin{subfigure}[t]{0.49\textwidth}
    \centering
    \includegraphics[width=\textwidth]{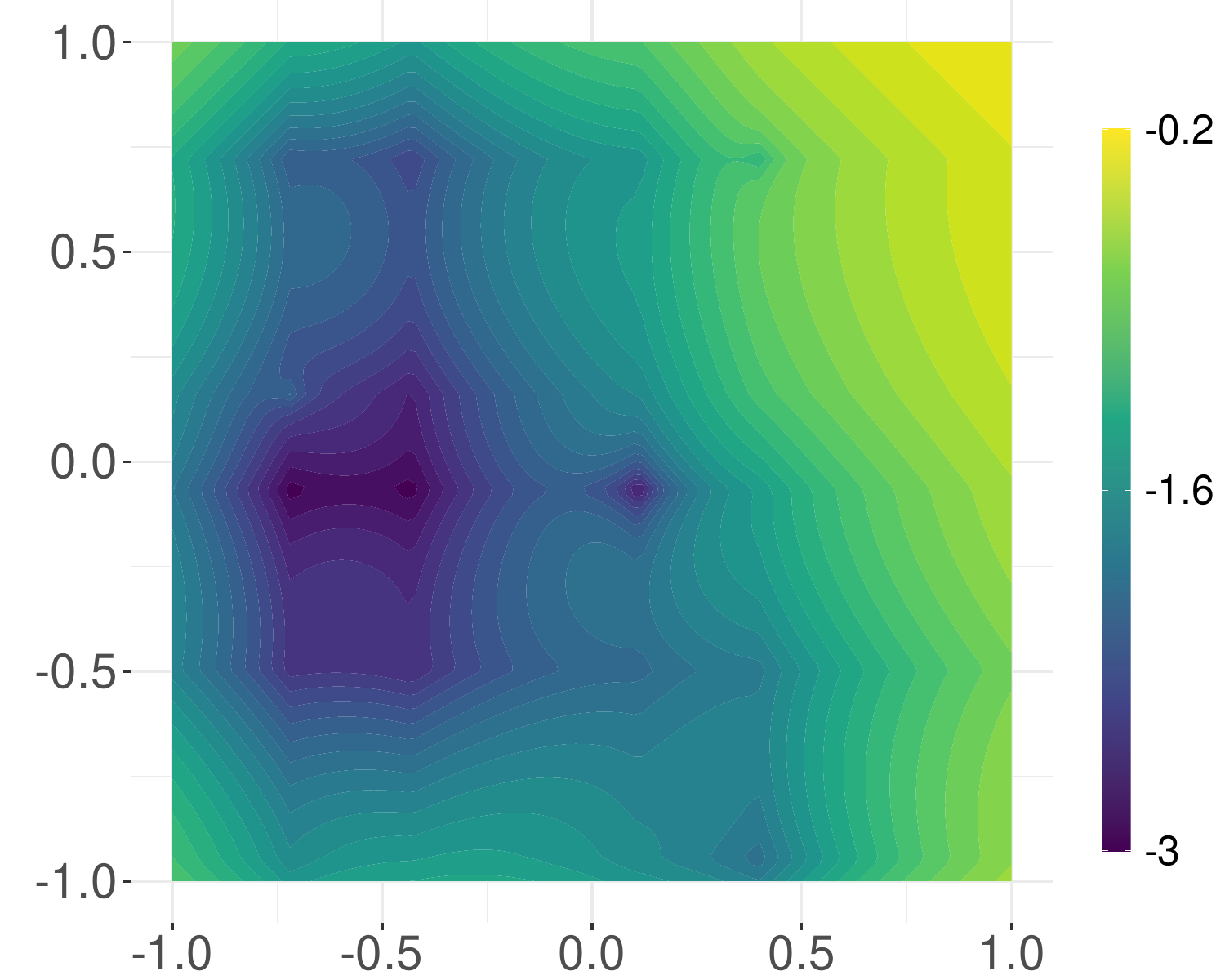}
    \caption{True function $f(\bfx)$ (contour).}
    \label{fig:matern52-true-contour}
  \end{subfigure}\hfill
  \begin{subfigure}[t]{0.49\textwidth}
    \centering
    \includegraphics[width=\textwidth]{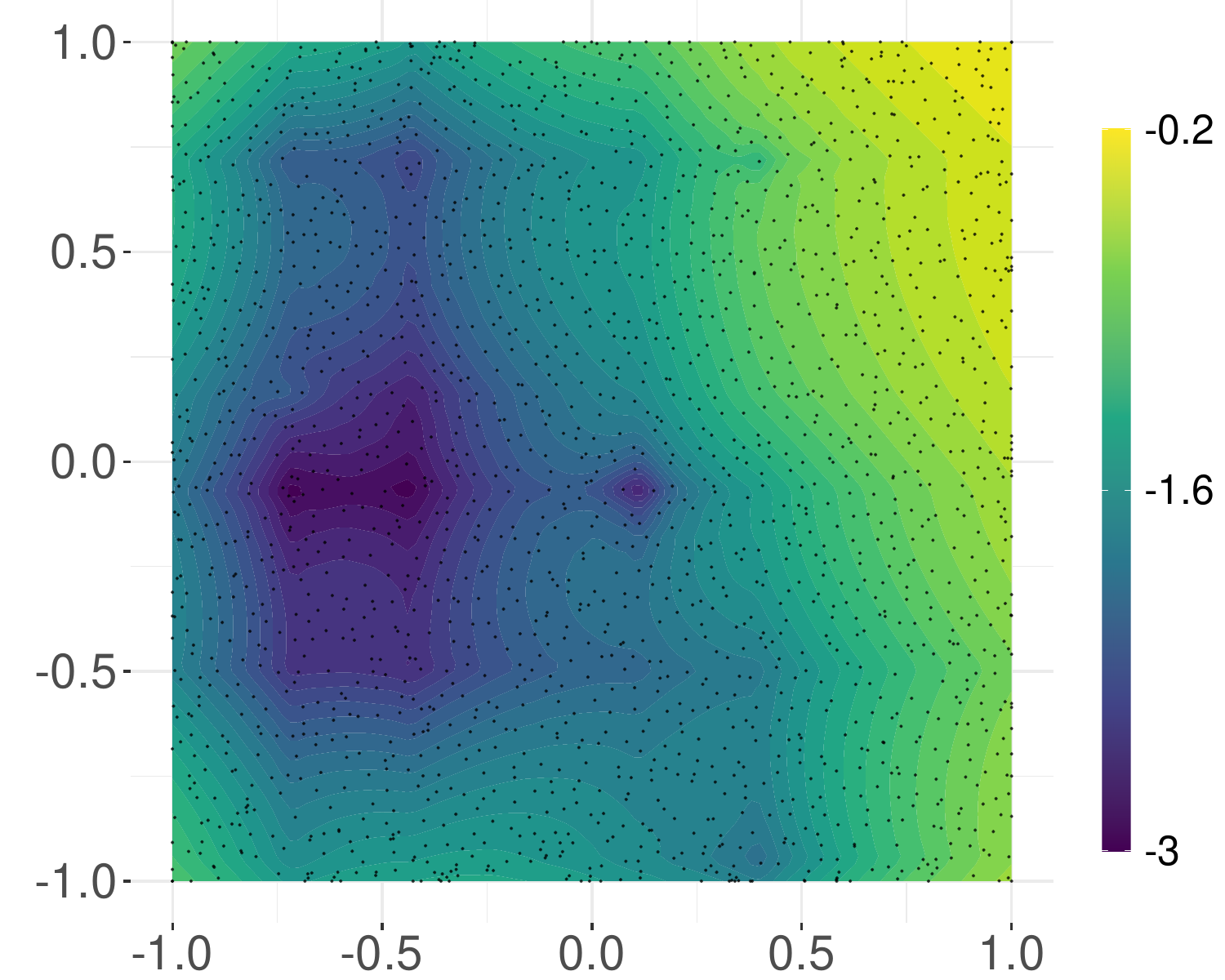}
    \caption{Posterior mean (contour).}
    \label{fig:matern52-mean-contour}
  \end{subfigure}
  \caption{Two-dimensional anisotropic Mat\'ern-$5/2$ example of HSGP-IMSE method.}
  \label{fig:prediction3d_matern5_2}
\end{figure}

\begin{figure}[!t]
  \centering
  \begin{subfigure}[t]{0.32\textwidth}
    \centering
    \includegraphics[width=\textwidth,]{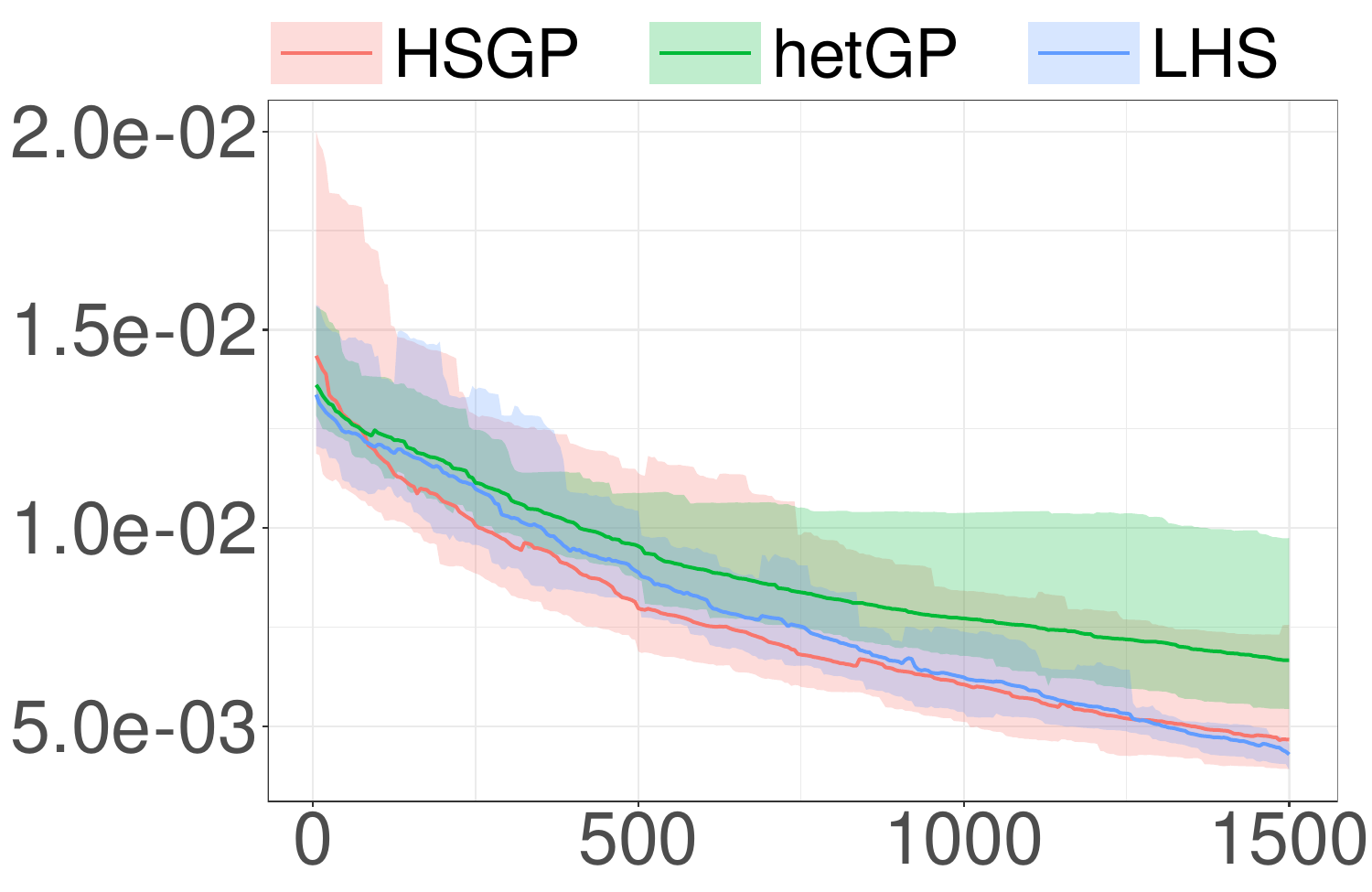}
    \caption{RMSE (Mat\'ern-$5/2$).}
    \label{fig:matern52-2d-rmse}
  \end{subfigure}\hfill
  \begin{subfigure}[t]{0.32\textwidth}
    \centering
    \includegraphics[width=\textwidth,]{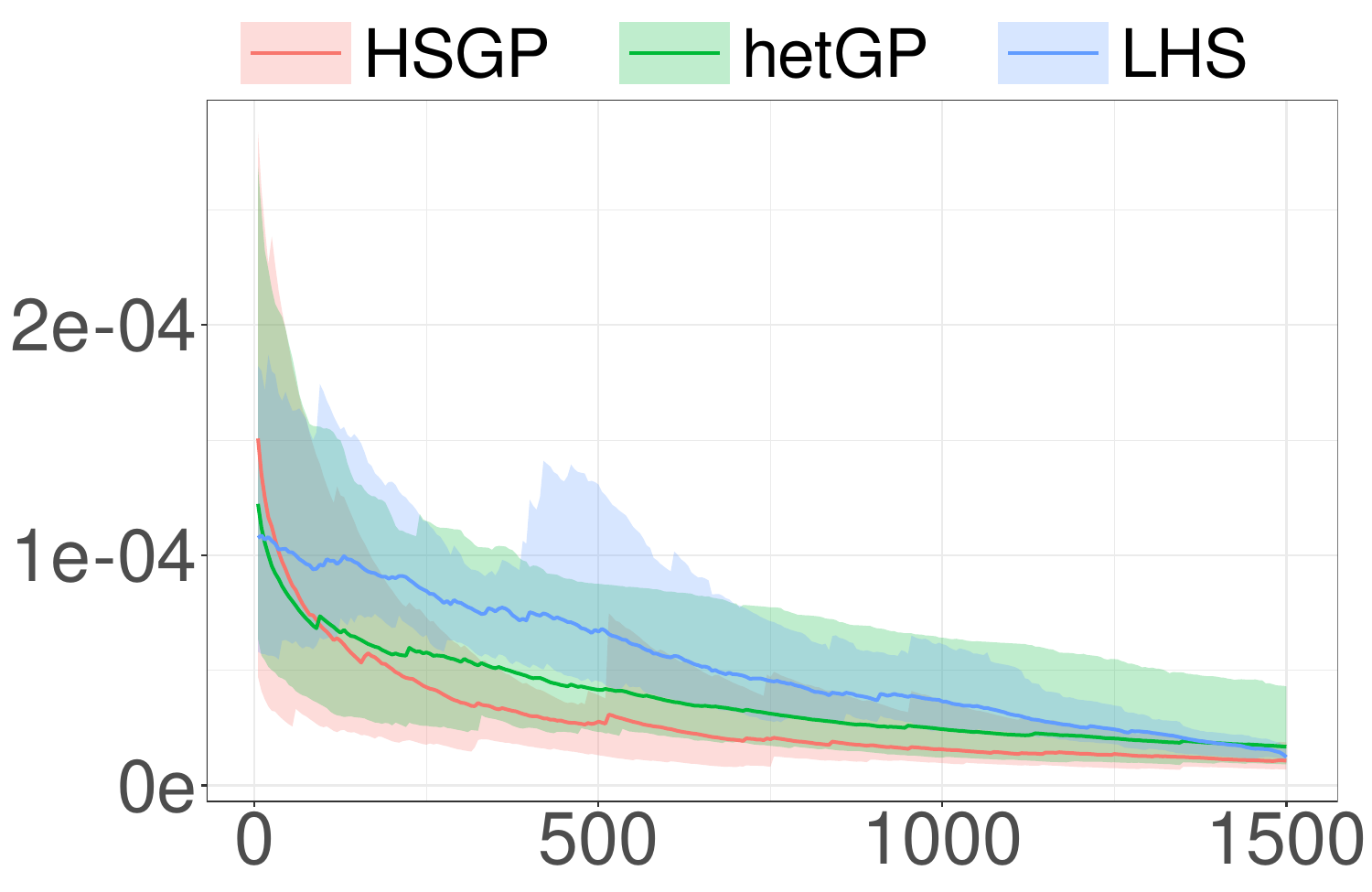}
    \caption{Variance (Mat\'ern-$5/2$).}
    \label{fig:matern52-2d-sd2}
  \end{subfigure}\hfill
  \begin{subfigure}[t]{0.32\textwidth}
    \centering
    \includegraphics[width=\textwidth,]{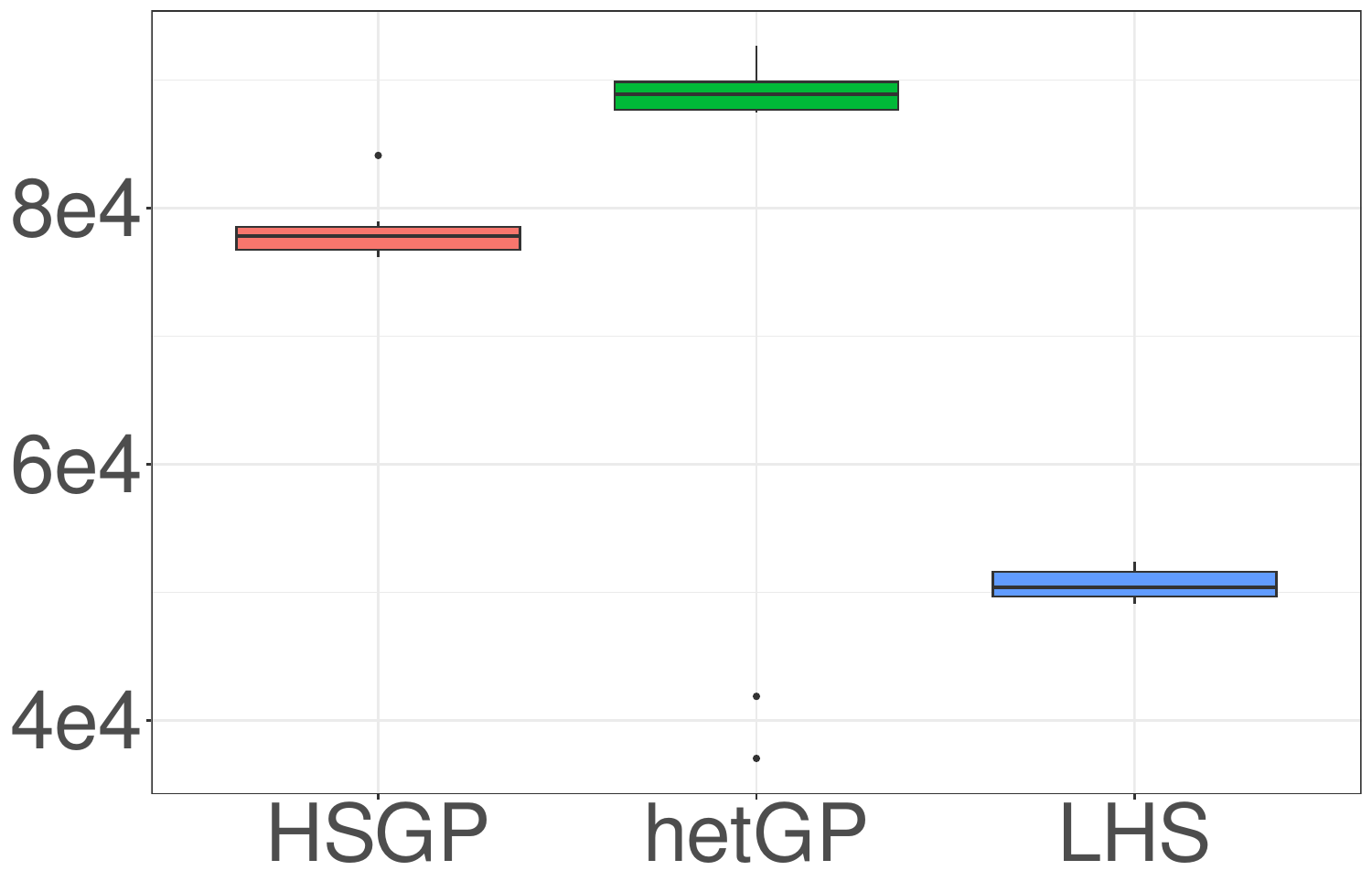}
    \caption{Time cost (Mat\'ern-$5/2$).}
    \label{fig:matern52-2d-time}
  \end{subfigure}
  \caption{Two-dimensional simulation with an anisotropic Mat\'ern-$5/2$ kernel.}
  \label{fig:sim_matern5_2_2d}
\end{figure}

Another test is training observations from $y = f_3(\bfx)$ using an anisotropic Mat\'ern-$5/2$ GP surrogate. The anisotropic Mat\'ern-$5/2$ kernel $k_A(\bfx, \bfx^\prime)$ is defined as 
\begin{align}
    k_A(\bfx, \bfx^\prime) = \prod_{k=1}^d k(x_k, x_k^\prime),
\end{align}
where $k(x_k, x_k^\prime)$ is defined in \eqref{eq:kernel_matern}. Figure~\ref{fig:prediction3d_matern5_2} and Figure~\ref{fig:sim_matern5_2_2d} report the results of fitting $f_3$ in $d=2$. The anisotropic setting stresses both model fitting and acquisition optimization: allowing dimension-specific length scales changes the posterior geometry, which can alter the location and sharpness of acquisition maxima. This test therefore verifies that $\widehat{\imse}_m$ tracks the same accuracy and uncertainty-reduction behavior as exact IMSPE while mitigating the increasing computational burden for GP surrogate with product kernels.

\begin{figure}[!t]
  \centering
  \begin{subfigure}[t]{0.32\textwidth}
    \centering
    \includegraphics[width=\textwidth,]{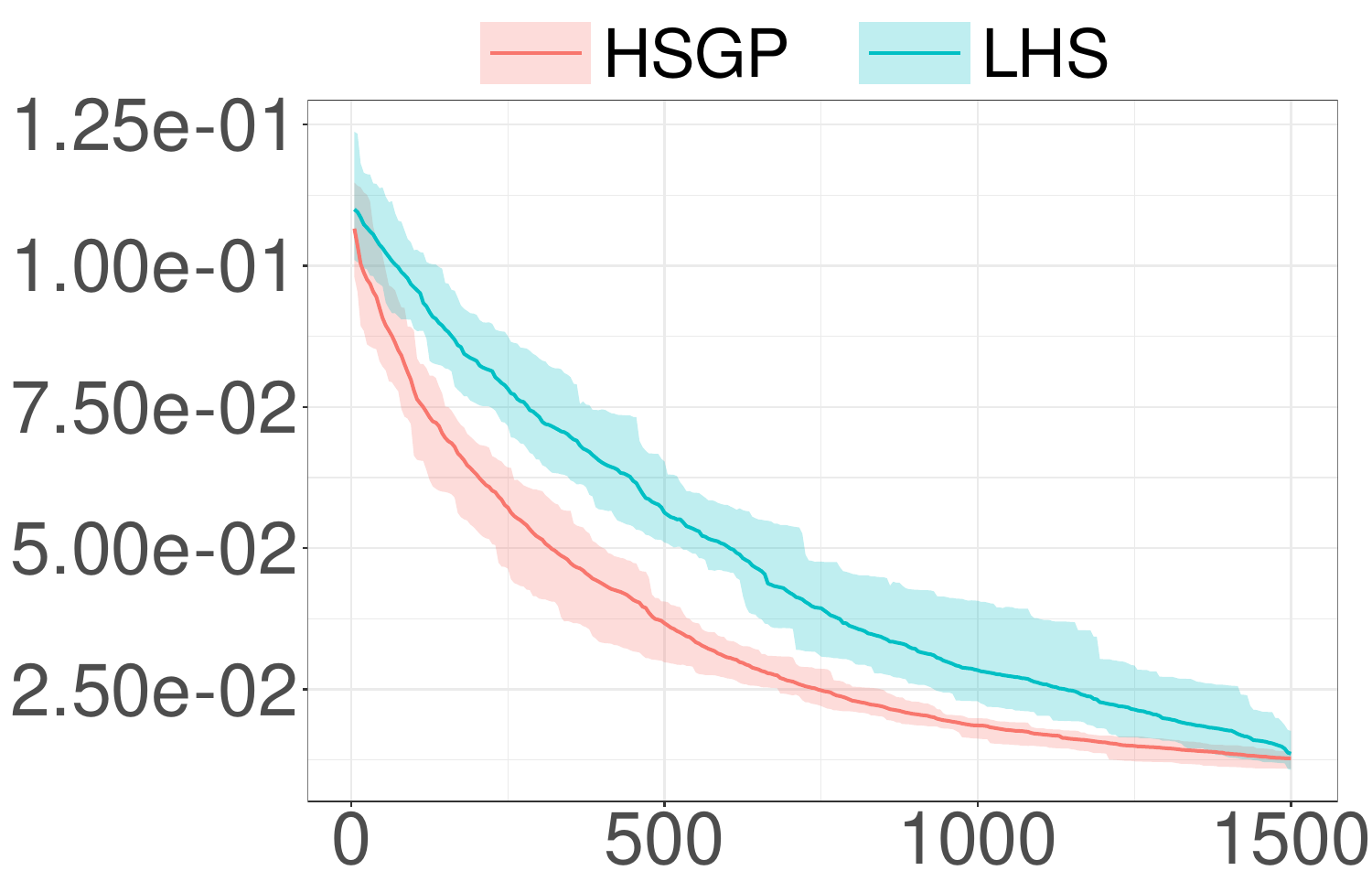}
    \caption{RMSE (Mat\'ern 2).}
    \label{fig:matern2-rmse}
  \end{subfigure}\hfill
  \begin{subfigure}[t]{0.32\textwidth}
    \centering
    \includegraphics[width=\textwidth,]{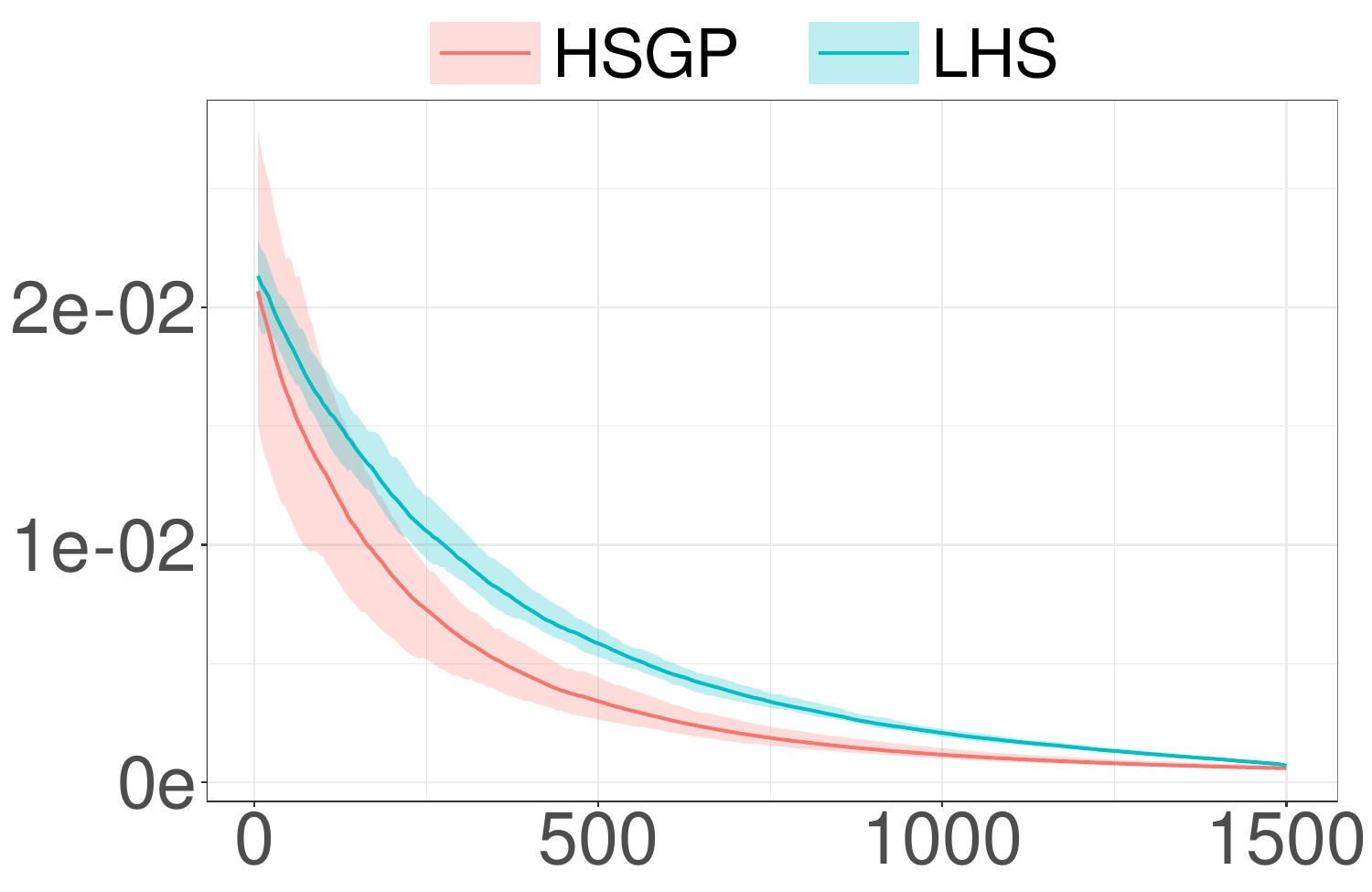}
    \caption{Variance (Mat\'ern 2).}
    \label{fig:matern2-sd2}
  \end{subfigure}
  \begin{subfigure}[t]{0.32\textwidth}
    \centering
    \includegraphics[width=\textwidth,]{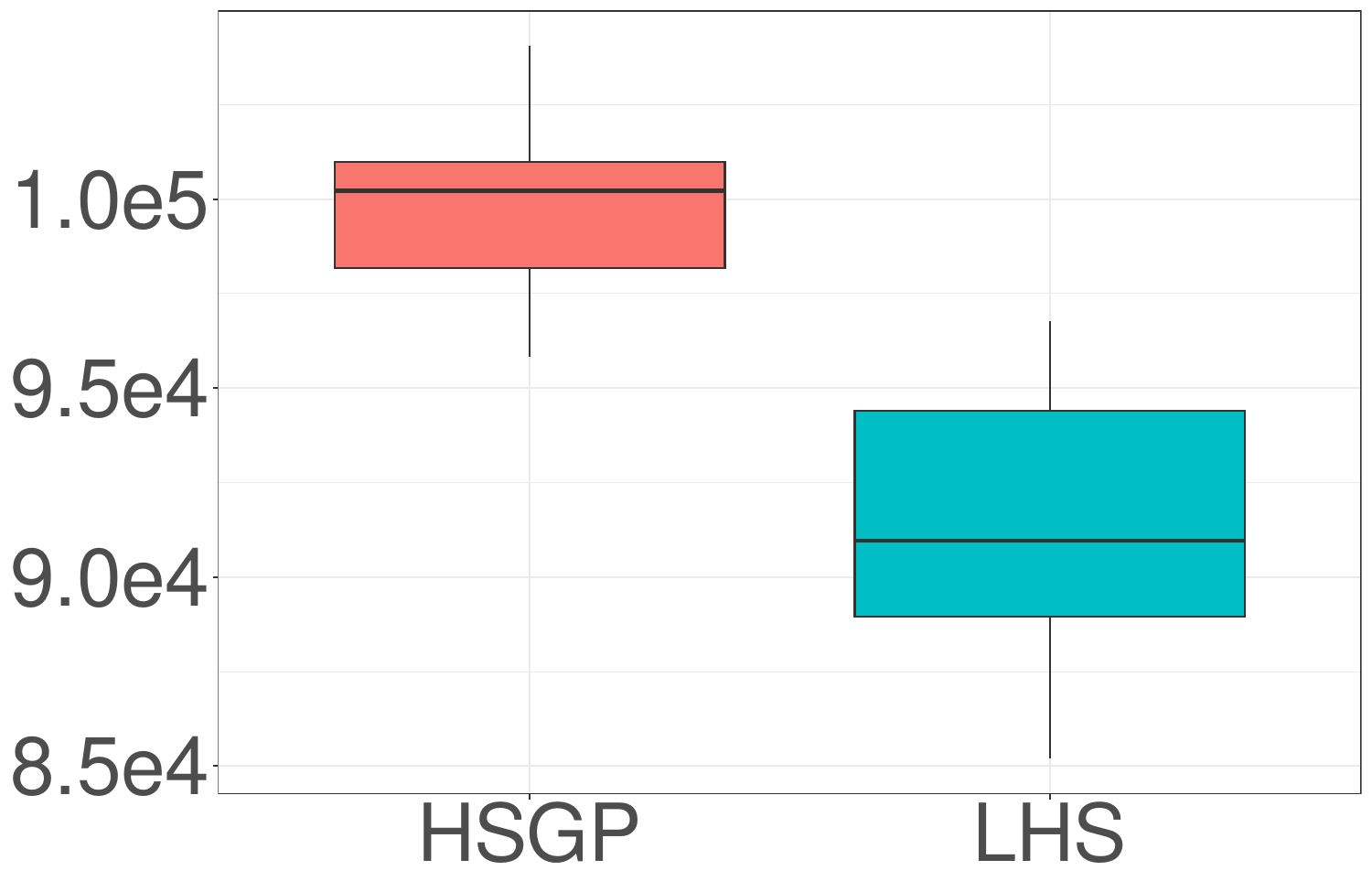}
    \caption{Time cost (Mat\'ern 2).}
    \label{fig:matern2-time}
  \end{subfigure}
  \caption{Two-dimensional simulation with an isotropic Mat\'ern-$2$ kernel.}
  \label{fig:sim_matern2_2d}
\end{figure}

\paragraph{Two-dimensional Mat\'ern-$2$ simulation.} Since the \texttt{hetGP} package only supports Gaussian, Mat\'ern-$3/2$, and Mat\'ern-$5/2$ kernels, the following experiments will only compare our HSGP-IMSE method with the LHS design. Figure~\ref{fig:sim_matern2_2d} reports the RMSE, mean posterior variance, and time cost for fitting $y = f_3(\bfx)$ samples using a Mat\'ern-$2$ kernel. The HSGP-IMSE achieves clearly lower RMSE and posterior variance with the same sample size, while its computational time is comparable to that of LHS.

\begin{figure}[!t]
  \centering
  \begin{subfigure}[t]{0.32\textwidth}
    \centering
    \includegraphics[width=\textwidth,]{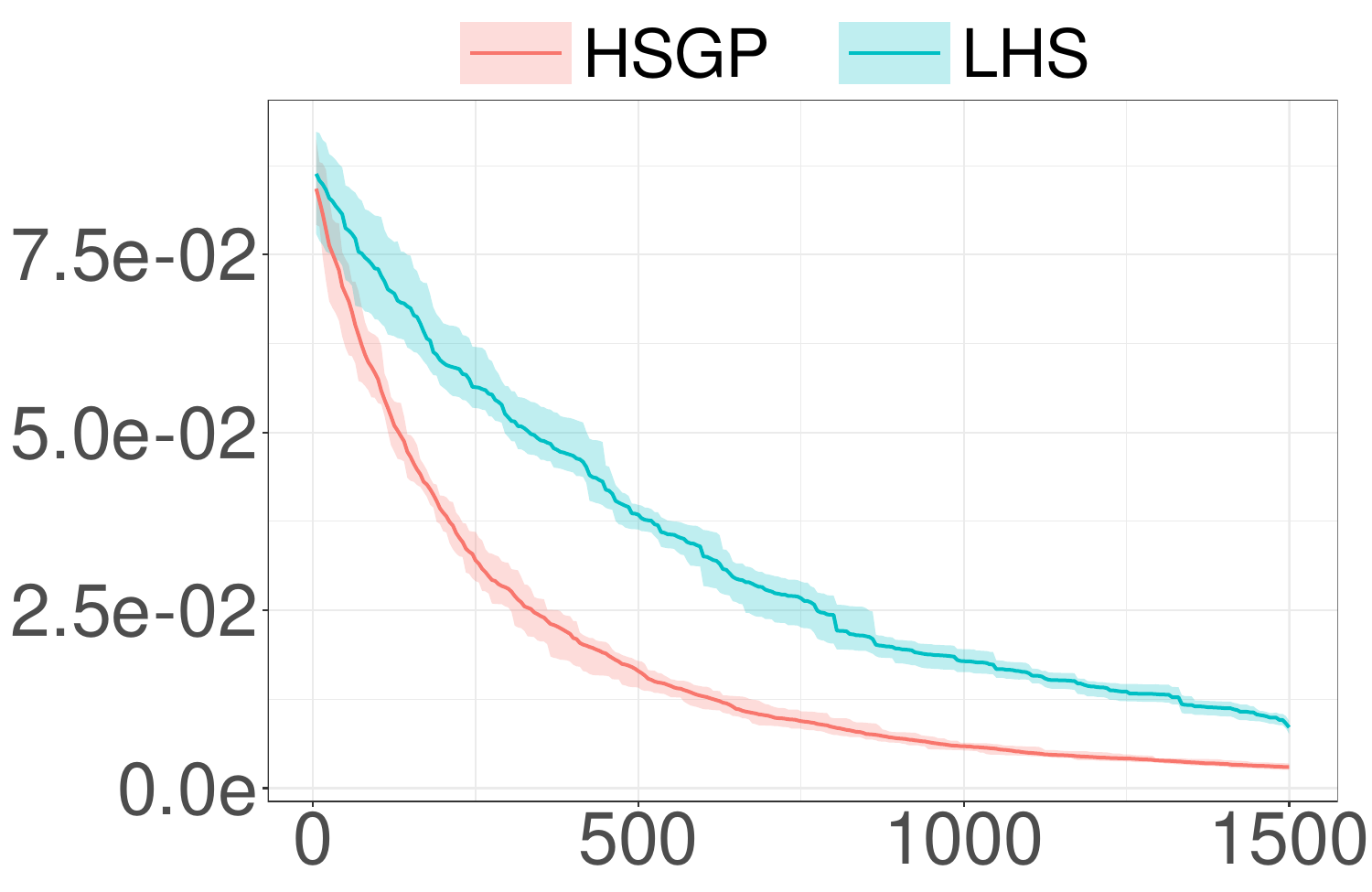}
    \caption{RMSE (GW).}
    \label{fig:gw-rmse}
  \end{subfigure}\hfill
  \begin{subfigure}[t]{0.32\textwidth}
    \centering
    \includegraphics[width=\textwidth,]{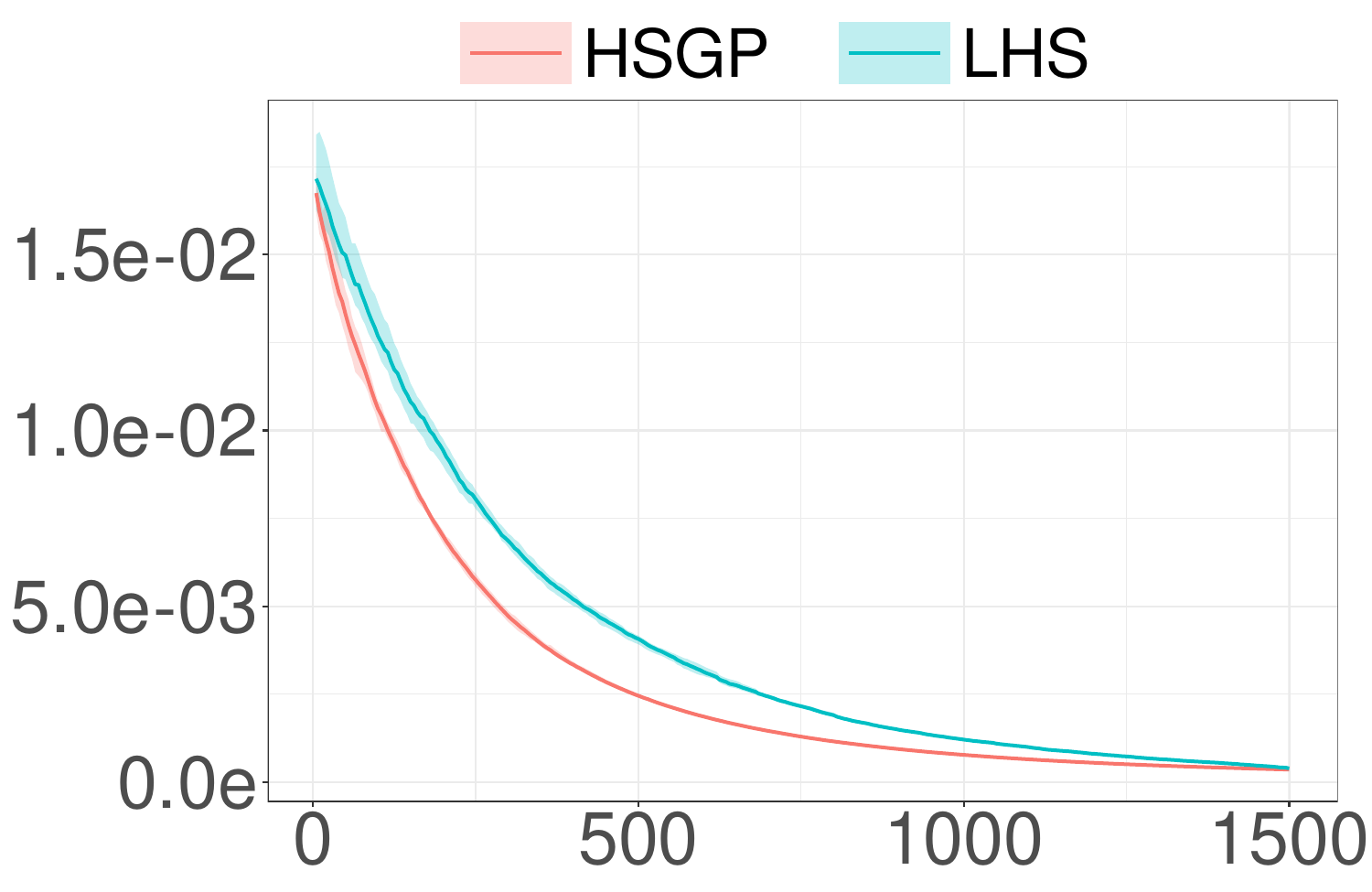}
    \caption{Variance (GW).}
    \label{fig:gw-sd2}
  \end{subfigure}
  \begin{subfigure}[t]{0.32\textwidth}
    \centering
    \includegraphics[width=\textwidth,]{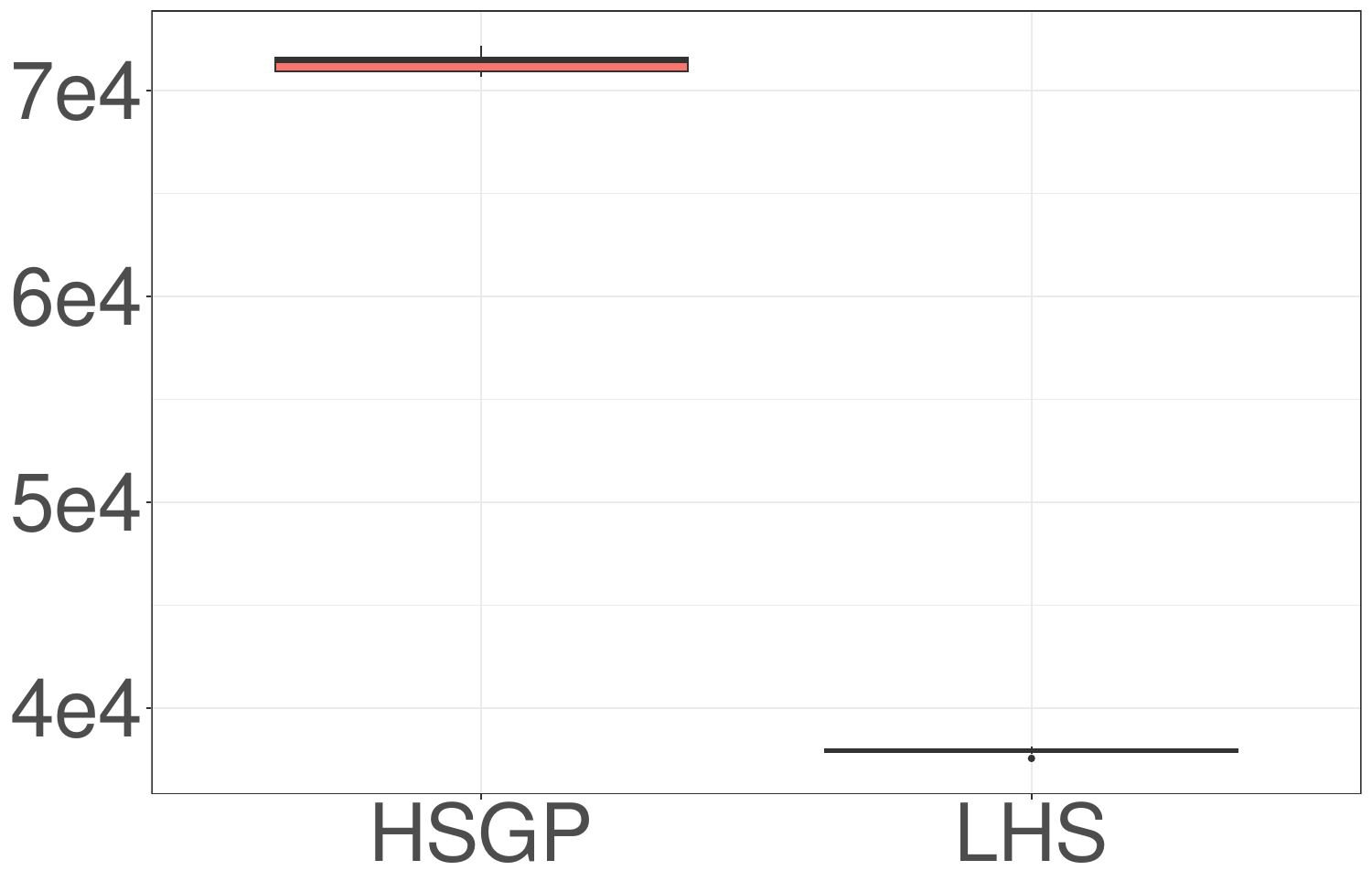}
    \caption{Time cost (GW).}
    \label{fig:gw-time}
  \end{subfigure}
  \caption{Two-dimensional simulation with a generalized Wendland (GW) kernel.}
  \label{fig:sim_gw_2d}
\end{figure}

\paragraph{Additional generalized Wendland kernel.} Finally, a generalized Wendland (GW) kernel (\citealt{bevilacqua:2019}) is considered to demonstrate that $\hsgp$-$\imse$ extends beyond the Gaussian and Mat\'ern families. Figure~\ref{fig:sim_gw_2d} reports the RMSE and mean posterior variance when fitting $f_4(\bfx)$ in $d=2$, under a GW covariance function with $\kappa=2$ and  $\mu=8.5$ in the parameterization of \citet{bevilacqua:2019}. According to \citet{bevilacqua:2019}, this GW kernel behaves similarly to the Mat\'ern-$5/2$ kernel. This experiment targets the methodological goal of generality: the same evaluation criteria (RMSE and uncertainty contraction, together with time cost) are used to verify that the surrogate-acquisition workflow remains effective when the kernel departs from standard kernel families such as Gaussian and Mat\'ern.

\section{Discussion} \label{sec:discussion}
This study proposes a fast and accurate HSGP-based approach for evaluating the $\imse$ acquisition in sequential design with GP surrogates. By replacing the original kernel with an HSGP approximation using parameter-free basis, the intractable integral admits a closed-form expression. This novel approach overcomes a key limitation of existing IMSE/IMSPE implementations, which are typically restricted to a narrow class of covariance families. In contrast, the proposed HSGP-IMSE applies broadly to kernels with closed-form spectral densities.

The work also provides a theoretical analysis of the approximation error. In particular, sharp exponential-decay bounds are established for the HSGP approximation of Gaussian and Mat\'ern kernels. Building on these results, convergence rates for the induced approximation of the IMSE acquisition with Mat\'ern kernels are derived in both noise-free and noisy settings.

Numerical experiments further demonstrate that HSGP-IMSE not only tracks the exact IMSE acquisition accurately, but also yields improved predictive performance in end-to-end sequential design, with lower RMSE, greater reduction in posterior uncertainty, and competitive computational cost. Overall, the proposed method extends the IMSE-based acquisition to a substantially broader class of GP surrogates, enabling flexible and computationally efficient sequential designs for global emulation.

Several directions remain for future work. While this study focuses on the Laplacian eigenbasis under Dirichlet boundary conditions, which yield closed-form trigonometric representations, alternative boundary conditions and associated bases in $L^2(\Omega,\PP)$ may also be feasible and could further mitigate boundary artifacts. In addition, the same techniques can be applied to problems involving integrals of kernel functions, such as kernel mean embedding, Bayesian quadrature, and spatial data fusion, which merit further investigation.

\vspace{6mm}

\noindent{\large \bf Acknowledgement}
\vspace{2mm}

This research was supported by the Singapore Ministry of Education Academic Research Funds Tier 1 Grant A-8002495-00-00.

\bibliographystyle{chicago}
\bibliography{ref}

\end{document}